\documentclass{amsart}

\usepackage{amsthm,amsfonts,amsmath,amssymb,latexsym,epsfig}
\usepackage{upref,amssymb,eucal,ae,enumitem}
\usepackage{mathdots}

\newtheorem{theorem}{Theorem}

\newenvironment{remark}{\medskip \refstepcounter{theorem}
\noindent  {\bf Remark \thetheorem}.\rm}{\,}

\def\<{\langle}
\def\>{\rangle}

\def\tm{\tilde{M}}
\def\tn{\tilde{n}}
\def\tg{\tilde{g}}
\def\mb#1{{\mathbb #1}}
\def\mc#1{{\mathcal #1}}
\def\mf#1{{\mathfrak #1}}
\def\BOne{{\mathchoice {\rm 1\mskip-4mu l} {\rm 1\mskip-4mu l}
                          {\rm 1\mskip-4.5mu l} {\rm 1\mskip-5mu l}}}

\begin{document}
\title[Conformal properties of spheres]
{Conformal properties of spheres}
\author{Santiago R. Simanca}
\email{srsimanca@gmail.com}

\begin{abstract}
We identify the cone of smooth metrics $\mc{M}(M)$ on a
closed manifold $M^n$ with the space of smooth isometric embeddings 
$f_g: (M,g) \rightarrow (\mb{S}^{\tn},\tg)$ into a standard sphere 
of large dimension $\tn=\tn(n)$, and their isotopic deformations, and the
space $\mc{C}(M)$ of conformal classes with the classes of 
metrics whose embeddings are isotopic to each other by conformal deformations.
Topological embeddings given by uniform 
limits of $C^1$ isometric embeddings of a metric, or isometric embeddings of 
metrics on the manifold with a different smooth structure, and their 
deformations, are carried by $(\mb{S}^{\tn},\tg)$ also, but when the latter
exist,  
they do not embed into a smooth flow of any $f_g(M)\hookrightarrow 
\mb{S}^{\tn}$. We characterize the metrics of constant scalar curvature in 
terms of properties of their isometric embeddings, and use this to prove a 
homotopy lifting property of the bundle $\mc{M}(M) \stackrel{\pi}{\rightarrow } 
\mc{C}(M)$ by Yamabe metrics. When $M$ carries an almost complex structure 
$J_0$, this result extends to a homotopy lifting property of the bundle 
$\mc{M}^{J_0}(M) \stackrel{\pi}{\rightarrow } \mc{C}^{J_0}(M)$ of 
metrics compatible with almost complex structures in the same 
orientation class as $J_0$, and their conformal classes, the lift now by almost 
Hermitian Yamabe metrics. We use these results and the gap theorem 
of Simons to study the existence and integrability properties of almost 
complex structures on spheres, and products. We 
find the sigma invariants of $Sp(2)$ (in fact, that of any compact Lie group) 
and the $M^7_k$ spheres of Milnor (in fact, that of any smooth $7$ dimensional
sphere), 
and with the exception of $\mb{P}^1(\mb{C})\times \mb{P}^1(\mb{C})$, 
the almost Hermitian sigma invariant of any product of spheres 
carrying almost complex structures. We organize the set of all manifolds 
with Yamabe, or almost Hermitian Yamabe metrics on them, in a Pascal  
like triangle that is set according to the symmetries of the metrics, and the 
values of their associated conformal invariants.
\end{abstract}
\maketitle

Our goal in this article is the study of almost complex structures on, 
and the computation of some of the sigma invariants of spheres and products.

In \S \ref{sec1}, we discuss the existence and integrability properties of 
almost complex structures on spheres, on the basis of technical 
results presented in \S \ref{sec2}. In \S \ref{sec2}, we characterize constant 
scalar curvature metrics on a manifold by properties of their Nash isometric 
embeddings into a standard sphere of sufficiently large dimension, and prove 
a homotopy lifting property by Yamabe metrics for the fibration of the space 
of metrics over the space of their conformal classes. When the manifold 
admits an  almost complex structure $J_0$, we characterize, likewise,  
metrics of constant almost Hermitian scalar curvature of metrics 
compatible with almost complex structures in the same orientation class 
as $J_0$, and extend the homotopy lifting property to the fibration of the 
space of these metrics over the space of their 
conformal classes, the lift now by almost Hermitian Yamabe 
metrics. We then carry out the discussion of \S \ref{sec1} for products of 
two spheres. In \S\ref{sec3}, we compute the sigma invariant of $Sp(2)$, and 
use its realizing metrics conveniently to find the sigma invariant of the 
$M^7_3$  sphere of Milnor, which implies that of all others 7 spheres, 
and continue by finding the almost Hermitian sigma invariant of 
all product of spheres that admit almost complex structures, with the 
exception of the K\"ahler surface $\mb{P}^1(\mb{C})\times \mb{P}^1(\mb{C})$
whose invariant we infer using only a qualitative argument. We end by  
organizing the information pertaining the invariants of all manifolds $M$ 
with Yamabe, and almost Hermitian Yamabe metrics on them, in a conformal 
Pascal triangle. The (known examples of) manifolds with metrics that realize 
their sigma invariants tie up its setting with the rational homotopy 
of the manifolds themselves. 
 
\section{Spheres with $\sqrt{-\BOne}$ tensors on their tangents} \label{sec1}
We denote by $(\mb{S}^n,g)$ the standard sphere in Euclidean 
$(\mb{R}^{n+1},\| \, \cdot \|^2)$ space, and by $\omega_n$ its volume;
when it becomes necessary, we shall refer to $g$ as $g_n$ also.
We distinguish notationwise a sphere of 
(optimal) dimension $\tn=\tn(n)$ that carries the Nash 
isometric embeddings of any $n$-Riemannian manifold $(M^n,g')$, and refer to it 
by $(\mb{S}^{\tn},\tg)$ instead. A metric $g'$ on $M$ is identified with its 
Nash isometric embedding $f_{g'}: (M,g') \rightarrow (\mb{S}^{\tn},\tg)$.
By the Palais' isotopic extension theorem, we can then identify 
the entire cone of metrics $\mc{M}(M)$ on $M$ with the space 
$\mc{M}_{\mb{S}^{\tn}}(M)$ of their Nash isometric embeddings 
into this fixed background, and their deformations. We use this device repeatedly to view Riemannian 
functionals 
in terms of functionals of the extrinsic quantities 
of the isometric embeddings of the metrics, 
and vice versa. 

Two major results, at opposite ends of this last interplay, single out the 
metric of $(\mb{S}^n,g)$ and its conformal class: The Yamabe problem on 
$(M^n,g')$ seeks metrics in the conformal class $[g']$ of $g'$, which 
achieve the infimum $\lambda(M,[g'])$ of the 
normalized total scalar curvature functional over the space of such
metrics, an invariant of the class $[g']$ that varies continuously with it.
We have the universal bound $\lambda(M,[g']) \leq \lambda(\mb{S}^n,[g]):=n(n-1)
\omega_n^{\frac{2}{n}}$, and the solution of the problem shows that 
$\lambda(M,[g'])$ is achieved by a
metric $g'\in [g']$ of constant scalar curvature $s_{g'}$, 
and that if $M \cong \mb{S}^n$ and $\lambda(M,[g'])=\lambda(\mb{S}^n,[g])$, 
then, $[g']=[g]$, and up to a conformal diffeomorphism, this solution is the 
standard metric $g$ on $\mb{S}^n$. On the other hand, for any $f_{g'} \in 
\mc{M}_{\mb{S}^{\tn}}(M)$ we have that $s_{g'}=n(n-1)+\| H_{f_{g'}}\|^2-\| 
\alpha_{f_{g'}}\|^2$, where 
$\alpha_{f_{g'}}$ and $H_{f_{g'}}$ are the second fundamental form and mean 
curvature vector of $f_{g'}$, respectively, and so the totally geodesic 
isometric embedding 
$f_g: (\mb{S}^n,g) \rightarrow (\mb{S}^{\tn},\tg)$ is minimal, and 
$s_g=n(n-1)$. By Simons' gap theorem, 
the space of immersions of $\mb{S}^n$ into $\mb{S}^{\tn}$
has a $C^2$ neighborhood of this $f_g \in \mc{M}_{\mb{S}^{\tn}}(\mb{S}^n)$ 
that contains no other minimal element.

We may bridge these two results, and see that if 
$t\rightarrow f_{g'_t} \in \mc{M}_{\mb{S}^{\tn}}(M)$ is
a path of isometric embeddings corresponding to a path of metrics  
$t\rightarrow g_t'\in \mc{M}(M)$, and that if 
$t\rightarrow f_{g'^{Y}_t}\in \mc{M}_{\mb{S}^{\tn}}(M)$ is an associated 
path of isometric embeddings of Yamabe metrics $g'^{Y}_t$ in 
$[g'_t]$, with $\mu_{g'^Y_t}(M)=\mu_{g'_t}(M)$, 
if $[g'^{Y}_{t}]=[g'_t]$ crosses the class $[g]$ at some point, then 
$M\cong \mb{S}^n$, and the entire path of classes 
$t\rightarrow [g'^{Y}_{t}]$ could had never been off the class 
$[g]$ before, or escape it after.

We give details of these, and various others defining terms and general 
technical results, in \S2. We use them now to prove the following.
\smallskip 

\begin{theorem} \label{th1}
{\rm (Borel \& Serre \cite{bs,bose}) (a) } 
No even dimensional sphere admits an almost complex structure unless its
dimension is $2$ or $6$.  {\rm (Poincar\'e \cite{poin}) (b)}
There exists a unique complex structure $J$ on $\mb{S}^2$, and as a complex 
manifold, $(\mb{S}^2,J) \cong \mb{P}^1(\mb{C})$. 
{\rm (\cite[Theorem 1]{sim5}) (c)} No almost complex 
structure on $\mb{S}^6$ is integrable.
\end{theorem}
\smallskip 

{\it Proof}. (a) Suppose that $J$ is an almost complex structure on the 
standard sphere $(\mb{S}^{n=2m},g)$, and consider the smooth path of metrics 
$$ 
[0,1]\ni t \rightarrow g_t( \, \cdot \, , \, \cdot \, ) = 
(1-t) g(J \, \cdot \, , J \, \cdot \, ) + 
t g( \, \cdot \, , \, \cdot \, ) \, .  
$$ 
We have that $g_1=g$, and so $[g_1]=[g]$. When $n=2$, $g_t=g$ and 
$[g_t]=[g]$ for all $t$. We prove that this last condition on the conformal
classes $[g_t]$ holds when $n>2$ also.

We consider the path 
$$
[0,1] \ni t \rightarrow f_{g_t}: ( \mb{S}^{n=2m},g_t)\rightarrow 
(\mb{S}^{\tn},\tg)
$$
of Nash isometric embeddings corresponding to the path of metrics 
$t\rightarrow g_t$. By Theorem \ref{thM2}, the path of 
conformal classes $t\rightarrow [g_t]$ admits
a lift to a path $t\rightarrow g_t^Y$ of Yamabe metrics   
$g_t^Y\in [g_t]$ such that that $\mu_{g_t^Y}=\mu_{g_t}$, with the
associated path of Nash isometric embeddings  
$$
[0,1] \ni t \rightarrow f_{g^Y_t}: ( \mb{S}^{n=2m},g^Y_t)\rightarrow 
(\mb{S}^{\tn},\tg)
$$
yielding paths of constant functions $t \rightarrow \| H_{f_{g_t^Y}}\|^2$ and 
$t \rightarrow \| \alpha_{f_{g_t^Y}}\|^2$ that are smooth and (at least) 
continuous, respectively. 
Each of the constant scalar curvature metrics $g_t^Y$ is either scalar flat, 
or Einstein, in which 
cases and if necessary, we can scale it by a constant $c>0$ to make 
$f_{c^2g_t^Y}$ minimal, or otherwise, its embedding $f_{g_t^Y}$ is minimal.
Irrespective of the situation, the indicated minimal isometric embedding
is a minimizer of the (intrinsic) functionals $\mc{W}_{f_{g}}(\mb{S}^{n})$ and 
$\mc{D}_{f_{g}}(\mb{S}^{n})$ in (\ref{eq24}) over the space of isometric 
embedding of metrics in the conformal class $[g_t]$. 
Since $g_1=g$, $g^Y_1=g$, and we have that 
$f_{g^Y_1}=f_g \in \mc{M}_{\mb{S}^{\tn}}(\mb{S}^n)$.

By Simons' gap theorem \cite[Theorem 5.2.3]{simo}, there exists
a $C^2$ neighborhood $\mc{U}_{f_g}$ of $f_g$ in the
space of immersions of $\mb{S}^n$ into $\mb{S}^{\tn}$ that contains
no other minimal element. Thus, $f_g$ is the only minimal element of 
$\mc{U}_{f_g} \cap \mc{M}_{\mb{S}^{\tn}}(\mb{S}^n)$. By the
solution to the Yamabe problem on the standard sphere, we then conclude that
for some $\varepsilon=\varepsilon(g)>0$, the restricted path
$[1-\varepsilon,1] \ni t \rightarrow f_{g^Y_t}$ consists of isometric
embeddings of metrics $g_t^Y$ that are smooth homothetic deformations of $g$, 
with the dilation factor determined by $\mu_{g_t}(\mb{S}^{n})$. 
These $g_t^Y$s are thus Einstein metrics in the standard conformal class 
$[g^Y_t]=[g_t]=[g]$.
Since the choice of $\varepsilon=\varepsilon(g)$ depends on (the conformal 
class of) $g$ only, by the continuity of $t\rightarrow [g_t]$ and that of 
$g\rightarrow \lambda(M,[g])$, we may rerun this argument starting with the
restriction of the path $[0,1]\ni t\rightarrow g_t^Y$ to, say, the interval 
$[0,1-0.9\, \varepsilon]$, and by iteration, prove that this result holds for 
all $t$s in the entire interval $[0,1]$. Thus, $[g^Y_t]=[g_t]=[g]$, as desired.

Since the path of classes $t \rightarrow [g_t]$ is constant, $[g_0]=[g]$, and 
$g_0$ and $g$ are conformally related. But $g_0= g\circ J$, so we can then
conclude that $g$ is $J$-invariant, and that the pair 
$(J,g)$ is an almost Hermitian structure on $\mb{S}^{n=2m}$. This implies that
$\mb{S}^{2m+1}$ is parallelizable \cite{kir}, and so $m=1$ or $m=3$.        

If we start now with an $n=2m$ sphere with an exotic differentiable structure, 
we pick an associated exotic diffeomorphism between this
sphere and $\mb{S}^{2m}$,
and by pull-back and push-forward, use it 
to transform tensor fields in the latter to tensor fields
in the former. In particular, any almost complex structure in one of these 
two spheres may be transformed into an almost complex structure in the 
other. The proof above then carries onto the exotic sphere, which can admit 
almost complex structures only if its dimension is either $2$ or $6$.   

(b) Given any metric $g'$ on $\mb{S}^2$, we denote by $J_{g'}$ the 
complex structure that it induces, so if we consider the path of
metrics
$$
[0,1]\ni t\rightarrow g_t=(1-t)g' + tg \, , 
$$
we obtain a smooth path of K\"ahler structures $t\rightarrow (J_{g_t},g_t)$ on 
$\mb{S}^2$ that starts at $(J_{g'},g')$ and ends at $(J_{g},g)$.  

We denote by $\omega_{g_t}$ the K\"ahler form of $(J_{g_t},g_t)$, and find
the minimizer of the squared $L^2$ form of the scalar curvature functional over
the space of K\"ahler metrics whose K\"ahler class is $[\omega_{g_t}]$, which
by the vanishing of the Kazdan-Warner obstruction \cite{kawa2},
must be of constant scalar curvature, and so an Einstein metric
in the conformal class $[g_t]$. We normalize it, if and when necessary, to an 
Einstein metric $g_t^n$ in $[g_t]$ such that 
$\mu_{g_t^n}(\mb{S}^{2})=\omega_2$, and
look at the ensuing path of Nash isometric embeddings    
$$
f_{g_t^n}: (\mb{S}^2,g_t^n)\hookrightarrow (\mb{S}^{\tn},
\tilde{g}) \, , 
$$
producing smooth paths of constant functions 
$t \rightarrow \| H_{f_{g_t^n}}\|^2$ and $t \rightarrow \| 
\alpha_{f_{g_t^n}}\|^2$. By the Gauss-Bonnet theorem, 
$$ 
s_{g_t^n}=\frac{8\pi}{\omega_2} =
2+ \| H_{f_{g_t^n}}\|^2 - \| \alpha_{f_{g_t^n}}\|^2 \, ,  
$$ 
and if we set $g_t^e=(1+\| H_{f_{g_t^n}}\|^2/4)g^n_t$, the family of 
isometric embeddings 
$$
f_{g_t^e}: (\mb{S}^2,(1+\| H_{f_{g_t^n}}\|^2/4)g_t^n )
\hookrightarrow (\mb{S}^{\tn},\tg)
$$
is minimal for all $t$, and the first of the associated paths 
of constant functions $t \rightarrow \| H_{f_{g_t}^e}\|^2$ and $t\rightarrow
\| \alpha_{f_{g_t^e}}\|^2$ is trivial. 

Since $g_1=g$, we have that $f_{g_1^n}=f_{g_1^e}=f_g$ and the continuous 
path $t\rightarrow \| \alpha_{f_{g_t^e}}\|^2$ vanishes at $t=1$. 
By the gap theorem of Simons \cite[Corollary 5.3.2]{simo}, this holds on
a nontrivial open neighborhood of $1\in [0,1]$, where we have that
$\| H_{f_{g_t^n}}\|^2=\| H_{f_{g_t^e}}\|^2=0=\| \alpha_{f_{g_t^e}}\|^2=
\| \alpha_{f_{g_t^n}}\|^2$, 
$[g_t^n]=[g_t^e]=[g_t]=[g]$, and modulo isometries, $f_{g_t^n}=f_{g_t^e}=f_g$.
If this neighborhood of $1$ were not equal to the entire interval $[0,1]$,
by Simons' gap theorem in its reinterpreted form as critical points of the 
squared $L^2$ of the mean curvature vector functional \cite[Theorem 2]{rss2},
we would then conclude that that for some $\bar{t}\in [0,1)$, the function 
$t\rightarrow \| \alpha_{f_{g_t^e}}\|^2$ would change discontinuously from
a value greater than $2$ on its immediate left, to $0$ on its 
immediate right.
Thus, $[g_t^n]=[g_t^e]=[g_t]=[g]$ for all $t$, and modulo 
isometries, $f_{g_t^n}=f_{g_t^e}=f_g$ always. Since
$g'=g_{t\mid_{t=0}}$, we have $[g']=[g]$ \cite[Theorem 1 a)]{sim2}, 
and as complex manifolds, $(\mb{S}^2,J_{g'})
=(\mb{S}^2,J_g) \cong \mb{P}^1(\mb{C})$ \cite[Corollary 2]{sim2}. 
\smallskip

(c) We let $(\mb{S}^6,J,g)$ be the standard octonionic almost Hermitian 
structure on $\mb{S}^6$. If $J'$ is any other almost complex structure on 
$\mb{S}^6$ in the same orientation class as that of $J$, we choose a path 
$[0,1]\ni t\rightarrow J_t$ of almost complex structures on $\mb{S}^6$ that 
begins at $J_0= J'$, and ends at $J_1 = J$, and define the associated 
path of $J_t$ invariant metrics 
$$
[0,1]\ni t \rightarrow g_{J_t}( \, \cdot \, , \, \cdot \, ) = \frac{1}{2}(
g( \, \cdot \, , \, \cdot \, )+ g( J_t \, \cdot \, , J_t \, \cdot \, ))\, . 
$$
We thus produce a path $t\rightarrow (J_t,g_{J_t})$ of almost Hermitian 
structures on $\mb{S}^6$ that begins at $(J',g_{J'})$ and ends at $(J,g)$. 
Since $g$ is $J$ invariant, we have that $g_{J_1}=g$, and $[g_{J_1}]=[g]$. 
The conformal class $[g_{J_0}]$ at the other end is a function of $J'$, and
since there are no orthogonal complex structures on $(\mb{S}^6,g)$ \cite{bl}, 
if $J'$ were to be integrable, we would have that $[g_{J_0}] \neq [g]$. We 
exclude this possibility by analyzing the path of conformal classes of the 
path of metrics   
$$
t\rightarrow g_t( \, \cdot \, , \, \cdot \, ) = 
(1-t) g_{J_0}( \, \cdot \, , \, \cdot \, ) + 
t g( \, \cdot \, , \, \cdot \, ) \, ,  
$$
for which we have that $[g_0]=[g_{J_0=J'}]$ and $[g_1]=[g]$, respectively, 
We adapt the argument in the early part of the proof of (a) above to
show that $[g_{t}]=[g]$ for all $t$. 

We consider the path of Nash isometric embeddings of the $g_t$s, 
$$
[0,1] \ni t \rightarrow f_{g_t}: ( \mb{S}^{n=6},g_t)\rightarrow 
(\mb{S}^{\tn},\tg) \, .
$$
By Theorem \ref{thM2}, 
the path $t\rightarrow [g_t]$ admits
a lift to a path $t\rightarrow g_t^Y$ of Yamabe metrics   
$g_t^Y\in [g_t]$ such that that $\mu_{g_t^Y}=\mu_{g_t}$, and the
associated path of Nash isometric embeddings  
$$
[0,1] \ni t \rightarrow f_{g^Y_t}: ( \mb{S}^{n},g^Y_t)\rightarrow 
(\mb{S}^{\tn},\tg)
$$
yields paths of constant functions $t \rightarrow \| H_{f_{g_t^Y}}\|^2$ and 
$t \rightarrow \| \alpha_{f_{g_t^Y}}\|^2$ that are smooth and (at least) 
continuous, respectively. By Theorem \ref{thM1}, the metrics $g_t^Y$ are 
either scalar flat, or Einstein, or otherwise, the embedding 
$f_{g_t^Y}$ is minimal, and irrespective of the case, either the isometric 
embedding itself or the isometric embedding of a suitable rescaling of the
metric, is a minimal minimizer of the functionals 
$\mc{W}_{f_{g}}(\mb{S}^{n})$ and $\mc{D}_{f_{g}}(\mb{S}^{n})$ in 
(\ref{eq24}) over the space of isometric embedding of metrics in the conformal 
class $[g_t]$. At $t=1$, $g^Y_1=g$, and $f_{g^Y_1}=f_g \in 
\mc{M}_{\mb{S}^{\tn}}(\mb{S}^n)$.

By Simons' gap theorem \cite[Theorem 5.2.3]{simo}, there exists
a $C^2$ neighborhood $\mc{U}_{f_g}$ of $f_g$ in the the space of immersions of 
$\mb{S}^{n=6}$ into $\mb{S}^{\tn}$ such that $f_g$ is the only minimal element 
of $\mc{U}_{f_g} \cap \mc{M}_{\mb{S}^{\tn}}(\mb{S}^n)$, and by the
solution to the Yamabe problem on the standard sphere, any Yamabe metric on
it coincides with $g$ up to a conformal diffeomorphism.
Hence, since either $f_{g_t^Y}$ or the isometric embedding of a suitable 
rescaling of 
$g_t^Y$ is minimal, there exists some $\varepsilon=\varepsilon(g)>0$ such that,
if $t\in [1-\varepsilon,1]$, the embedding $f_{g_t^Y}$ is either $f_{g}$ itself,
or it is the embedding of a metric $g_t^Y$ that is homothetic to
$g$, with the scale factor determined by $\mu_{g_t}(\mb{S}^6)$.
Either way, these $g_t^Y$s are Einstein, and $[g_t^Y]=[g_t]=[g]$.
Since the choice of $\varepsilon$ depends on $g$ only, and the functions
$t\rightarrow [g_t]$ and $g \rightarrow \lambda(M,[g])$ are continuous, 
we can rerun this argument for the restriction of the path of 
metrics to say, the interval 
$[0,1-0.9\, \varepsilon]$, and by iteration, prove that these results hold
for all $t$ in $[0,1]$. Thus, the entire path 
$t\rightarrow g_t^Y$s is a path of homothetics of $g$, hence a path of
Einstein metrics homothetic to $g$, and $[g_t^Y]=[g_t]=[g]$ for all 
$t\in [0,1]$. It follows that $[g_0]=[g_{J_0=J'}]=[g]$, 
contradicting their earlier indicated inequality.
\qed
\smallskip

\section{Metrics of constant scalar curvature: Yamabe metrics}
\label{sec2}
For any (closed) manifold $M=M^n$, the space $\mc{M}(M)$ of smooth 
Riemannian metrics on $M$ is the space of sections of the bundle of positive 
definite symmetric two tensors over $M$, provided with the induced compact 
open smooth topology on the sections of the bundle of all symmetric tensors 
\cite{ebin}. The space $\mc{C}(M)$ of conformal classes of metrics 
is the quotient space by conformal equivalences of metrics, which we
provide with the induced quotient compact open smooth topology. 
If $\mc{M}_{[g]}(M)$ denotes 
the set of metrics in the conformal class $[g]$, the quotient fibration
\begin{equation} \label{ma}
\mc{M}(M) \stackrel{\pi}{\rightarrow} \mc{C}(M) 
\end{equation}
has fiber $\mc{M}_{[g]}(M)$ over 
$[g]\in \mc{C}(M)$, and we have the foliated decomposition
$$
\mc{M}(M)=\cup_{[g]\in \mc{C}(M)}\mc{M}_{[g]}(M) \, . 
$$

The spaces $\mc{M}(M)$ and $\mc{C}(M)$ are contractible, each of them being a 
locally trivial Serre fibration over $M$ with contractible fiber, and 
(\ref{ma}) itself is a Serre fibration with a special kind of homotopy lifting 
property, which follows by the existence of special metrics of constant 
scalar curvature in each conformal class. 

Indeed, by the Nash isometric embedding theorem \cite{nash}, any 
Riemannian $n$-manifold $(M^n,g)$ can be isometrically embedded into 
a standard sphere $(\mb{S}^{\tn}, \tg) \hookrightarrow  
(\mb{R}^{\tn+1}, \| \phantom{.} \|^2)$ in Euclidean space
of sufficiently large, but fixed, optimal dimension $\tn= \tn(n)$. This
$\tn$ depends on $n$ but not on $M$.
  
For any path   
\begin{equation} \label{emc}
[0,1]\ni t \rightarrow f_{g_t}: (M,g_t) \mapsto (\mb{S}^{\tn}, \tg)
\end{equation}
of Nash isometric embedding deformations of $f_g=f_{g_t}\mid_{t=0}$ 
corresponding to a
path $t\rightarrow g_t$ of metric deformations of $g$, by the Palais' isotopic
extension theorem \cite{pal}, there exists a smooth path 
of diffeomorphisms
$$
F_t: \mb{S}^{\tn} \rightarrow \mb{S}^{\tn}
$$
such that $F_t(f_g(x))=f_{g_t}(x)$, and since the metrics on the submanifolds 
are all induced by the metric $\tg$ on the fixed background $\mb{S}^{\tn}$, we 
obtain by restriction a diffeomorphism 
\begin{equation} \label{emc2}
F_t\mid_{f_g(M)}: (f_g(M),\tg) \rightarrow (f_{g_t}(M), \tg) \, .   
\end{equation}
Thus, we may identify the open cone $\mc{M}(M)$ of smooth Riemannian 
metrics on $M$ with the space of their isometric embeddings into 
$(\mb{S}^{\tn},\tg)$, 
\begin{equation} \label{met}
\mc{M}_{\mb{S}^{\tn}}(M)=\{f_{g'}: \; f_{g'}: (M,g') 
\rightarrow (\mb{S}^{\tn},\tg) \; \text{isometric embedding}\, \}\, ,   
\end{equation}
and their isotopic deformations. The metric $g'$ in $\mc{M}(M)$ corresponds 
to the end point of the
path of isometric embeddings $F_t\mid_{f_g(M)}=f_{g_t}$ connecting it to 
$f_g$, modulo isometries of the background sphere, and so $g$ and $g'$ 
are associated to the submanifolds $f_g(M)$ and $f_{g'}(M)$ of $\mb{S}^{\tn}$ 
ambiently isotopic to each other. A submanifold $f_g(M)$ 
that is ambiently isotopic to itself is associated to the same metric $g$, in 
which case the said isotopy reflects the right action by pull-back of the 
identity component of the diffeomorphism group of $M$, realized by 
restrictions of ambient space 
diffeomorphisms, on the metric. The fixing of just one $f_g(M)$ 
fixes the differentiable structure of all the other submanifolds $f_{g'}(M)$,  
ambiently isotopic to it, that we consider in this way. 

The topological embeddings of the manifold $M$ that are uniform limits of 
$C^1$ isometric embeddings of a smooth metric, or ``exotic'' 
embeddings of smooth metrics on the manifold provided with a different 
smooth structure, are also carried by $(\mb{S}^{\tn},\tg)$, and 
coexist with the $f_g$s above, and their deformations. Since any 
neighborhood $U$ of the identity diffeomorphism $\BOne_{f_g(M)}$ contains a 
neighborhood $V$ such that no element of $V$ embeds in a topological flow
\cite{pali,pali1}, the deformations of these topological or exotic isometric 
embeddings, when they exist, do not embed into a smooth flow of the 
submanifolds, and are in this sense, far away from the deformations of any 
$f_g$ in $\mc{M}_{\mb{S}^{\tn}}(M)$.

Suppose now that (\ref{emc}) is a path of embeddings 
deformations of $f_g$ corresponding to a path
\begin{equation} \label{cd}
[0,1]\ni t \rightarrow g_t =e^{2\psi (t)}g
\end{equation}
of conformally related metrics on $M$, where $\psi(0)=0$. Then the metric 
$g$ on $f_g(M)$ and $g_t$ on $f_{g_t}(M)$ are globally conformally related to
each other, and since  
the tensor $F_t^* \tg$ is just the background metric $\tg$ on
$\mb{S}^{\tn}$ acted on by the diffeomorphism $F_t$, and 
the metrics on the submanifolds are conformally related, there exists
a function $u(t)$ such that
$$
F_t^* (\tg\mid_{f_{g_t}(M)}) = e^{2 u(t)(f_g(\, \cdot \,))} \tg \mid_{f_g(M)} =
e^{2 u(t)(f_g(\, \cdot \,))} \tg\mid_{f_{g}(M)} \, ,
$$
where the factor $e^{2 u(t)}$ and that in (\ref{cd}) are related 
to each other by $e^{2\psi(t)(\, \cdot \, )}=e^{2u(t)\circ f_g(\, \cdot \, )}$. 
If $u(t)$ is conveniently extended to a tubular neighborhood of $f_g(M)$, 
we may then view (\ref{emc}) as the family of pointwise conformal isometric 
embeddings 
\begin{equation}\label{embt}
f_{g_t}: (M,e^{2u(t)\circ f_g}g) \rightarrow (f_g(M),e^{2u(t)(f_g(\, \cdot
\, ))}\tg)
\hookrightarrow (\mb{S}^{\tn}, \tg) 
\end{equation}
of the fixed submanifold $f_g(M)$ of $\mb{S}^{\tn}$ with a varying conformal
metric on it. The metric $g'=g_t\mid_{t=1}$ in the conformal class of $[g]$ 
corresponds 
to the end point of this conformal embedding deformation of $f_g$, modulo 
isometries of the background sphere.   

In order to elaborate on the geometry of the embeddings of the conformally 
related metrics (\ref{embt}), and their relations to $f_g$, we recall briefly 
some material from \cite[\S 2]{gracie} and  
\cite[\S 3]{sim5} that we shall need for doing so.
In general, if $(\tm^{\tn},\tg)$ is some fixed Riemannian 
background, and 
\begin{equation} \label{emb}
f_g : (M,g)\rightarrow (\tm,\tg)
\end{equation}
is an isometric embedding of $(M^n,g)$ into $(\tm^{\tn},\tg)$,   
Gauss' equation reads 
\begin{equation}
g(R^g(X,Y)Z,W)= \tg(R^{\tg}(X,Y)Z,W)+\tg(\alpha(X,W),\alpha(Y,Z))-
\tg(\alpha(X,Z),\alpha(Y,W))\, , 
\label{cu}
\end{equation}
where $R^g$ and $R^{\tg}$ are the curvature tensors, and $\alpha:=\alpha_{f_g}$ 
is the second fundamental form of $f_g(M)$. If
$\{ e_1, \ldots ,e_n, \nu_{n+1}, \ldots, \nu_{\tn}\}$ is a (local) 
orthonormal frame of $\tm$, with $\{ e_1, \ldots, e_n\}$ an orthonormal frame 
of $M$, the Ricci tensors of $g$ and $\tg$ are thus related by
\begin{equation}
\begin{array}{rcl} \label{ri}
r_g(X,Y) & = & 
r_{\tg}(X,Y) - \sum_{i=n+1}^{\tn} \tg(R^{\tg}(\nu_i,X)Y,\nu_i) +
\tg(H,\alpha(X,Y))-
\\ & & \mbox{} \hspace{2in} \sum_{i=1}^n \tg(\alpha(e_i,X),\alpha(e_i,Y)) \, ,
\end{array}
\end{equation}
where $H:=H_{f_g}=\sum_i\alpha(e_i,e_i)$ is the mean curvature vector, 
and the scalar curvatures by  
\begin{equation} \label{sc}
\begin{array}{rcl}
s_g & = & \sum_{i,j\leq n}K^{\tg}(e_i,e_j)
+\tg(H,H)- \tg(\alpha,\alpha) \, , 
\end{array}
\end{equation}
where $\sum_{i,j}K^{\tg}(e_i,e_j)$ is the extrinsic scalar curvature
on $f_{g}(M)$.

By (\ref{sc}), the total scalar curvature decomposes as
$$
\begin{array}{rcl}
{\displaystyle \int s_g d\mu_g } & = & {\displaystyle  
\Theta_{f_g}(M):=\int \sum_{i,j\leq n}K_{\tg}(e_i,e_j)d\mu_g  
+ \Psi_{f_g}(M):=\int \| H_{f_g}\|^2 d\mu_g } \\ & &  
- {\displaystyle \Pi_{f_g}(M):=\int \| \alpha_{f_g}\|^2 d\mu_g }\, ,
\end{array}
$$
and if $N=2n/(n-2)$, the Yamabe functional is the linear 
combination
\begin{equation} \label{func} 
\lambda_g(M):=\frac{1}{\mu_g^{\frac{2}{N}}}\int s_g d\mu_g =
\frac{1}{\mu_g^{\frac{2}{N}}}\Theta_{f_g}(M)+
\frac{1}{\mu_g^{\frac{2}{N}}}\Psi_{f_g}(M) -
\frac{1}{\mu_g^{\frac{2}{N}}}\Pi_{f_g}(M) \, .  
\end{equation}
The critical points of this $\lambda_g$ in the conformal class  
are the metrics 
of constant scalar curvature in the class, and there is always at least one
that realizes the intrinsic conformal invariant
$$ 
\lambda(M,[g])= \inf_{g'\in \mc{M}_{[g]}(M)} \lambda_{g'}(M)=
\inf_{g'\in \mc{M}_{[g]}(M)} \frac{1}{\mu_{g'}(M)^{\frac{2}{N}}}\int s_{g'} 
d\mu_{g'} \, ,
$$ 
\cite{ya,au,tr,sc}. These metric minimizers are by definition the 
Yamabe metric representatives of $[g]$. We have the universal bound \cite{au}
$$
\lambda(M,[g])\leq \lambda(\mb{S}^n,[g_n])=
\lambda_{g_n}(\mb{S}^n)=n(n-1) \omega_n^{\frac{2}{n}} \, ,
$$
and so 
\begin{equation} \label{sig}
\sigma(M)= \sup_{[g]\in \mc{C}(M)} \lambda(M,[g]) 
\end{equation}
is a well-defined differentiable invariant of $M$ \cite{sc2}. 

We may compute relevant intrinsic and extrinsic quantities of 
the varying conformal deformation $f_{g_t}$ in (\ref{embt}) in terms of 
likewise quantities for $f_g$ itself, and suitable local operators acting on 
$u(t)$. If for convenience we set $\tilde{g}_t=e^{2u(t)}\tilde{g}$,
we obtain that 
\begin{equation}  \label{pti}
\begin{array}{rcl} 
\sum_{i,j} K^{\tilde{g}_t}(e_i^t, e_j^t) & = &  e^{-2u(t)}\left(
\sum_{i,j} K^{\tilde{g}}(e_i,e_j) -2(n-1)\left({\rm div}_{f_g(M)}\nabla^{\tg} 
u^{\tau}- \right. \right. \\& & \left. \left. \hspace{2cm}\tg(H_{f_g},
\nabla^{\tg} u^\nu)-\tg(du^{\tau},du^{\tau})+\frac{n}{2}\tg(du,du) 
\right) \right) \, , \vspace{1mm} \\
\| H_{f_{g_t}}\|^2 & = & e^{-2u(t)}(\| H_{f_g}\|^2 
- 2n\tg(H_{f_g},\nabla^{\tg}u ^\nu)
+n^2 \tg(\nabla^{\tg}u ^\nu, \nabla^{\tg} u^\nu)) \, , \vspace{1mm} \\
\| \alpha_{f_{g_t}}\|^2 & = & e^{-2u(t)}(\| \alpha_{f_g}\|^2 - 2\tg(H_{f_g},
\nabla^{\tg}u^\nu) +n \tg(\nabla^{\tg}u ^\nu, \nabla^{\tg} u^\nu)) \, ,
\end{array}
\end{equation}
where to emphasize the role of the constant sectional curvature of the 
background sphere, we have kept in place the notation 
$\sum_{i,j} K^{\tg_t}(e_i^t,e_j^t)$ instead of its actual value $n(n-1)$, and
the superindices $\tau$ and $\nu$ refer to the tangential and normal components.
By (\ref{sc}), these relations imply the intrinsic scalar curvature equation
\begin{equation} \label{gra} 
s_{g_t}= e^{-2u(t)}\left( s_{g}-2(n-1){\rm div}_{f_g(M)}(\nabla^{\tg}u)^\tau
-(n-1)(n-2)\tg(\nabla^{\tg}u^{\tau},\nabla^{\tg}u^\tau)  \right) \, .
\end{equation}
With these identities in place, we may characterize metrics of constant scalar 
curvature in terms of properties of their associated isometric embeddings 
(\ref{met}), and among them, the Yamabe metrics in particular.

Along volume preserving deformations of any kind,  
$\left( \frac{1}{\mu^{\frac{2}{N}}_{g_t}}
\Theta_{f_{g_t}}(M) \right)$ remains constant, and so its variation is
trivial. If the deformations are in addition conformal, 
we use the top identity in (\ref{pti}) to compute
this variation in terms of extrinsic quantities of $f_g$, and variations of 
the function $u(t)$, and find it to be 
$$
{\displaystyle \frac{n-2}{
\mu^{\frac{2}{N}}_{g_t}}\int_{f_g(M)}\hspace{-0.7cm} e^{(n-2)u}\dot{u}\left( 
n(n-1)\! - \!\frac{1}{\mu_{g_t}}\int_{f_{g_t}} \hspace{-2mm}
n(n-1)d\mu_{g_t}\right)d\mu_{\tg} \! - \!  
\frac{2(n-1)}{\mu^{\frac{2}{N}}_{g_t}}\int_{f_g(M)} \hspace{-0.7cm} 
e^{(n-2)u(t)}\tilde{g}(n\nabla^{\tilde{g}}u^\nu -H_{f_g},
\nabla^{\tilde{g}} \dot{u}^\nu ) d\mu_{\tilde{g}}} \, ,   
$$
thus obtaining that the weak equation
\begin{equation}\label{basic}
0=\frac{d}{dt} \left( \frac{1}{\mu^{\frac{2}{N}}_{g_t}}
\Theta_{f_{g_t}}(M) \right) =
\frac{-2(n-1)}{\mu^{\frac{2}{N}}_{g_t}}\int_{f_g(M)} 
e^{(n-2)u(t)}\tilde{g}(n\nabla^{\tilde{g}}u^\nu -H_{f_g},
\nabla^{\tilde{g}} \dot{u}^\nu ) d\mu_{\tilde{g}} 
\end{equation} 
holds for all $t$. 

Proceeding similarly using the last two identities in (\ref{pti}) instead,
we find that 
\begin{equation} \label{gr32}
\begin{array}{rcl}
{\displaystyle \frac{d}{dt} \left(\frac{1}{\mu^{\frac{2}{N}}_{g_t}}
{\Psi}_{f_{g_t}}(M) \right)} \hspace{-2mm} & = & \hspace{-2mm}
{\displaystyle  \frac{n-2}{
\mu^{\frac{2}{N}}_{g_t}}\int \dot{u}\left(
\| H_{f_{{g}_t}}\|^2-
\frac{{\Psi}_{f_{g_t}(M)}}{\mu_{g_t}}\right)
d\mu_{\tilde{g}_t} 
+
\frac{2n}{ \mu^{\frac{2}{N}}_{g_t}}
\int e^{(n-2)u(t)}\tilde{g}(n\nabla^{\tilde{g}}u^\nu \! -H_{f_g},
\nabla^{\tg} \dot{u}^\nu ) d\mu_{\tilde{g}}
\, ,} \\ 
{\displaystyle \frac{d}{dt} \left(\frac{1}{\mu^{\frac{2}{N}}_{g_t}}
{\Pi}_{f_{g_t}}(M)\right)} \hspace{-2mm}  & = & \hspace{-2mm}   
{\displaystyle  \frac{n-2}{
\mu^{\frac{2}{N}}_{g_t}}\int \dot{u}\left(
\| \alpha_{f_{{g}_t}}\|^2-
\frac{{\Pi}_{f_{g_t}(M)}}{ \mu_{{g}_t}}\right)
d\mu_{\tilde{g}_t} 
+ 
\frac{2}{ \mu^{\frac{2}{N}}_{g_t}}
\int e^{(n-2)u(t)}\tilde{g}(n\nabla^{\tilde{g}}u^\nu \! -H_{f_g},
\nabla^{\tg} \dot{u}^\nu ) d\mu_{\tilde{g}}
\, .} 
\end{array}
\end{equation}
By the previous result, the last summand on the right side of these 
expressions vanishes for all $t$. Hence, along volume preserving conformal 
deformations $f_{g_t}$ of $f_g$, these functionals have a critical point at 
time $t$ if, and only if, the first summand on the right side is equal to 
zero for all possible variations $\dot{u}$ of $u$ at this $t$, and when this 
happens, we must have that both  
$\| H_{f_{g_t}}\|^2$ and $\| \alpha_{f_{g_t}}\|^2$ are constant 
functions at the said $t$.
These are precisely the critical points of the Yamabe functional
(\ref{func}), the metrics of constant scalar curvature in $[g]$.

Since $(\mb{S}^{\tn},\tg)$ carries all the isometric embeddings of metrics
on $M$, and $\tg$ has constant sectional curvature one,  
for any $f_g \in \mc{M}_{\mb{S}^{\tn}}(M)$, 
the function $\| H_{f_g}\|^2$ is constant \cite[Theorem 6]{sim5}\footnotemark, 
and therefore, along volume 
\footnotetext{This fact is true for any $C^2$ isometric embedding that is a
deformation of smooth volume preserving isometric embeddings in the conformal
class. But there may be other isometric embeddings of the same metric with
nonconstant function $\| H_{f_g}\|^2$, even for $g$s of constant 
scalar curvature (in which case, $\| H_{f_g}\|^2-\| \alpha_{f_g}\|^2$ would
be constant). Such an embedding cannot be $C^2$ close to, hence 
deformable into one $f_g\in \mc{M}_{\mb{S}^{\tn}}(M)$, preserving the volume 
and conformal class.}      
preserving conformal deformations, $f_{g_t}$ is a critical point of 
$\frac{1}{\mu_{g_t}^{\frac{2}{N}}}\Phi_{f_{g_t}}(M)$ for all $t$.    
Hence, the metric $g_t\mid_{t=1}$ is a critical point of the Yamabe 
functional in its conformal class if, and only if,  $\| \alpha_{
f_{g_t\mid_{t=1}}}\|^2$ 
is a constant function \cite[Theorem 7]{sim5}. 
Notice that though the codimension of the constant
scalar curvature embedding $f_{g_t\mid_{t=1}}$ may be strictly less that 
$\tn-n$, $f_{g_t\mid_{t=1}}(M)$ is volume preserving deformable to any smooth 
submanifold of the form $f_{g'}(M)$, $g'\in \mc{M}_{[g]}(M)$ of the same
volume, and in fact, it is homeomorphically isotopic to any 
$f_{g'}(M)$ for any $C^1$ embedding $f_{g'}$ of such a metric $g'$, or  
their topological uniform limit deformations. Thus, the characterization given 
of constant scalar curvature metrics applies only to the embeddings in 
$\mc{M}_{\mb{S}^{\tn}}(M)$, and not, for instance, to a smooth embedding 
in low codimension that is not a volume preserving deformations of any 
other embedding $f_{g'}\in \mc{M}_{\mb{S}^{\tn}}(M)$, $C^2$ close to it, and
in the same conformal class.

If $s_g$ is a positive constant, and the embedding $f_g$ is minimal, 
then $g$ is a Yamabe metric \cite[Theorem 4]{sim5}, while for $(M,g)$s such
that $\lambda(M,[g])\leq 0$, if $s_g$ is a nonpositive constant, 
$g$ is, up to homothetics, the unique Yamabe metric in its class.
We sharpen these results next.

None of the three functionals on the right side of (\ref{func}) is
intrinsically defined on the space of isometric embeddings of metrics in 
$\mc{M}_{[g]}(M)$, but we may rearrange their linear combination into two
functionals instead, and resolve this issue. For let us consider the $n$ 
dimensional versions of the 
extrinsic functionals of embedded surfaces used in \cite{sim2}, 
\begin{equation} 
\begin{array}{rcl}
\mc{W}_{f_{g}}(M) & = & {\displaystyle 
\frac{n}{n-1}\Theta_{f_{g}}(M)+\Psi_{f_{g}}(M) 
}\, ,\vspace{1mm}\\
\mc{D}_{f_{g}}(M) & = & {\displaystyle 
\frac{1}{n-1}\Theta_{f_{g}}(M)+\Pi_{f_{g}}(M) 
}\, ,  
\label{eq11}
\end{array}
\end{equation}
and observe that under conformal deformations $f_{g_t}$ of $f_g$,  
we have that
\begin{equation} \label{eq23}
\begin{array}{rcl}
\mc{W}_{f_{g_t}}(M) \! \! \! \!  &=& \! \! \! \! {\displaystyle 
\int_{f_g(M)} \hspace{-0.55cm}  e^{(n-2)u} \! 
\left( \frac{n}{n-1}\sum_{i,j}K^{\tilde{g}}(e_i,e_j)+ \| H_{f_g}\|^2 
+n (n-2)\| du^{\tau}\|^2 \right) d\mu_{g} }\, , \vspace{1mm} \\
\mc{D}_{f_{g_t}}(M) \! \! \! \! &= & \! \! \! \! { \displaystyle 
  \int_{f_g(M)} \hspace{-0.55cm} e^{(n-2)u} \! 
\left( \frac{1}{n-1}\sum_{i,j}K^{\tilde{g}}(e_i,e_j)+ \| \alpha_{f_g}\|^2 
+(n-2)\| du^{\tau}\|^2 \right) d\mu_{g} }\, ,   
\end{array}
\end{equation}
the densities of which are local operators of the function $u\mid_{f_{g}(M)}$ 
only, in contrast with what happens with their defining summands, which depend 
on $\nabla^{\tg}u^{\nu}$ also. We obtain two intrinsically well-defined 
functionals 
\begin{equation} \label{eq24}
\begin{array}{rcl}
 & &  \frac{1}{\mu_{g}^{\frac{2}{N}}}\mc{W}_{f_g}(M) \\ 
 & \nearrow & \\
\mc{M}_{[g]}(M) \ni g &  &  \\
 & \searrow & \\
 & &  \frac{1}{\mu_{g}^{\frac{2}{N}}}\mc{D}_{f_g}(M) 
\end{array}
\end{equation}
that are invariant under conformal transformations of the ambient 
background $(\mb{S}^{\tn},\tg)$, and (\ref{func}) restricted to 
$\mc{M}_{[g]}(M)$ decomposes now as 
$$ 
\lambda_g(M)= \frac{1}{\mu_{g}^{\frac{2}{N}}}\mc{W}_{f_g}(M)-
\frac{1}{\mu_{g}^{\frac{2}{N}}}
\mc{D}_{f_g}(M) \, .  
$$ 
(The analogous decomposition when $n=2$ is 
$4\pi \chi(M)=\mc{W}_{f_g}(M)-\mc{D}_{f_g}(M)$, with 
$\mc{W}_{f_g}$ and $\mc{D}_{f_g}$ each constant in the 
conformal class, and so of trivial variation. We then find 
a canonical representative of each class by finding critical points
of the squared $L^2$ norm of $s_g$, and fixing their areas conveniently,  
find the minimizers of $\mc{W}_{f_g}$ and $\mc{D}_{f_g}$ across 
classes \cite[Theorems 1 \& 9]{sim2} (cf. proof of Theorem 1 (b) here).)

\begin{theorem} \label{thM1}
Let $(M^{n},g)$ be a Riemannian manifold. A metric $g'\in [g]$ 
has constant scalar curvature if, and only if,
its isometric embedding $f_{g'}$ is a critical point of 
$\frac{1}{\mu_{g}^{\frac{2}{N}}}\mc{W}_{f_g}(M)$ and 
$\frac{1}{\mu_{g}^{\frac{2}{N}}}\mc{D}_{f_g}(M)$ in the
space of conformal deformations of $f_g$. $g'$ is a Yamabe metric if, and only 
if, $f_{g'}$ is a minimizer of these functionals in the said space of
deformations, which happens, if, and only if, $g'$ is either Einstein,
or scalar flat, in which cases there exists $c>0$ such that 
$f_{c^2g'}$ is minimal
and a minimizer of $\mc{W}_{f_g}(M)$ and $\mc{D}_{f_g}(M)$ in the conformal 
class, or if otherwise, $g'$ has constant scalar curvature and 
$f_{g'}$ is minimal and a  
minimizer of $\mc{W}_{f_g}(M)$ and $\mc{D}_{f_g}(M)$ in the conformal 
class.
\end{theorem}
\smallskip
 
{\it Proof}. If $g'$ has constant scalar curvature, 
then $f_{g'}$ is a critical point of the functionals $\Psi_{f_g}$ and 
$\Pi_{f_g}$ in the space of volume preserving conformal deformations of the 
embedding, of constant density each \cite[Theorems 6 \& 7]{sim5}. Since 
the linear combinations (\ref{eq11}) of these functionals define 
$\mc{W}_{f_g}(M)$ and $\mc{D}_{f_g}(M)$, respectively, $f_{g'}$ must be a 
critical point of these two functionals also in the said space of conformal
deformations, and therefore, a critical point of the volume normalized 
functionals in the space of conformal deformations of $f_g$, volume preserving 
or not. 

In order to prove the converse, we first use (\ref{eq23}) to 
write the variations of the volume normalized versions of $\mc{W}$ and $\mc{D}$
in terms of the variation $\dot{u}$ of the path $u=u(t)$ associated to
the conformal deformation. We obtain
$$ 
\begin{array}{rcl}
{\displaystyle \frac{d}{dt} \frac{1}{\mu^{\frac{2}{N}}_{\tilde{g}_t}}
\mc{W}_{f_{\tilde{g}_t}}(M)} & = & 
{\displaystyle  \frac{n-2}{\mu^{\frac{2}{N}}_{\tilde{g}_t}} \left(
\int_{f_g(M)} \hspace{-0.2cm} e^{nu}\dot{u}\left( e^{-2u}\left( \frac{n}{n-1} 
\sum K^{\tilde{g}}(e_i,e_j) + \| H_{f_{g}}\|^2
\right.  \right. \right. } \vspace{1mm}\\ & & 
-{\displaystyle \left.  \left. \left. 
n(n-2)\tg(du^{\tau}, du^{\tau}) + 2n  
\Delta^g u \right)  - 
\frac{1}{\mu_{\tilde{g}_t}}
\mc{W}_{f_{\tilde{g}_t}}(M)\right)  d\mu_{g} \right)   
 }\, , \vspace{1mm} \\
{\displaystyle \frac{d}{dt} \frac{1}{\mu^{\frac{2}{N}}_{\tilde{g}_t}}
\mc{D}_{f_{\tilde{g}_t}}(M)} \! \! 
& = & \! \!  {\displaystyle  \frac{n-2}{\mu^{\frac{2}{N}}_{\tilde{g}_t}}
\int_{f_g(M)} \hspace{-0.2cm} e^{nu}\dot{u}\left( e^{-2u}\left( \frac{1}{n-1} 
\sum K^{\tilde{g}}
(e_i,e_j) + \| \alpha_{f_{g}}\|^2
\right. \right. } \\ & & 
- {\displaystyle \left. \left. \left. 
(n-2)\tg(du^{\tau},du^{\tau}) + 2 \Delta^{g} u \right) 
-\frac{1}{\mu_{\tilde{g}_t}}\mc{D}_{f_{\tilde{g}_t}}(M)
  \right) d\mu_{g}  \right) }    
\, ,  
\end{array}
$$ 
respectively. It follows that if a volume preserving conformal deformation 
$f_{g_t}$ is a critical point of $\mc{W}_{f_g}$ and $\mc{D}_{f_g}$ at some 
$t$, there must be constants $\lambda_{\mc{W}}$, $\lambda_{\mc{D}}$ such that 
$$
\begin{array}{rcl}
e^{-2u}\left( n^2 + \| H_{f_{g}}\|^2 \! -n(n-2)\tg(du^{\tau},du^{\tau})
+ 2n \Delta^{g} u \right) - {\displaystyle \frac{1}{\mu_{\tilde{g}_t}}
\mc{W}_{f_{\tilde{g}_t}}}(M) \! \! \! & = & \! \! \! 
{\displaystyle \lambda_{\mc{W}}} \, = \, 0 \, ,   
\vspace{1mm} \\
e^{-2u}\left( n + \| \alpha_{f_{g}}\|^2 \! -(n-2)\tg(du^{\tau},du^{\tau})
+ 2 \Delta^{g} u \right) - {\displaystyle \frac{1}{\mu_{\tilde{g}_t}}
\mc{D}_{f_{\tilde{g}_t}}(M)} \! \! \! & = & \! \! \! 
{\displaystyle \lambda_{\mc{D}}}\, = \, 0 \, , 
\end{array}
$$    
hold at the said $t$, the indicated values of the constants being 
obtained by integrating the resulting equation against the measure 
$e^{nu(t)}d\mu_g$. 
Thus, this $f_{g_t}$ is a critical point of the functionals
${\displaystyle \frac{1}{\mu_{g}^{\frac{2}{N}}}\mc{W}_{
f_{g}}(M)}$ and  ${\displaystyle \frac{1}{\mu_{g}^{\frac{2}{N}}}
\mc{D}_{f_{g}}(M)}$ in the
space of conformal deformations, volume preserving or not.   
As a linear combination of $\Theta_{f_g}(M)$ and $\Psi_{f_g}(M)$, the
functional $\mc{W}_{f_g}(M)$ is stationary along volume preserving conformal 
deformations $f_{g_t}$. 
Thus, if $f_{g_t}$ is a critical point of $\mc{D}_{f_g}(M)$ along volume
preserving conformal deformations at some $t$, this $f_{g_t}$ is a also 
a critical point of $\mc{W}_{f_g}(M)$ in the said space, and therefore,
a critical point of $\frac{1}{\mu_g^{2/N}}\mc{D}_{f_g}(M)$ and  
$\frac{1}{\mu_g^{2/N}}\mc{W}_{f_g}(M)$ in the space of conformal deformations.

Let us say that $f_{g_t}\mid_{t=1}$ is such a critical point, and for
convenience, let us set $g'=\tilde{g}_t\mid_{t=1}=e^{2u(1)}\tg$. By
(\ref{gra}), if we take the difference of the critical point equations above 
at $t=1$, we obtain that
$$
\begin{array}{rcl}
s_{g'} 
& = & 
{\displaystyle e^{-2u(1)}\left( s_g  
+ 2(n-1) \Delta^{g} u -
(n-1)(n-2)\tg(du^{\tau},du^{\tau}) \right)\mid_{t=1} } \\ 
& = & 
{\displaystyle e^{-2u(1)}\left( n(n-1) + \| H_{f_{g}}\|^2 - 
\| \alpha_{f_{g}}\|^2+ 2(n-1) \Delta^{g} u -
(n-1)(n-2)\tg(du^{\tau},du^{\tau}) \right)\mid_{t=1} } \\ 
& = & 
{\displaystyle \frac{1}{\mu_{g'}} \mc{W}_{f_{g'}}(M)
- \frac{1}{\mu_{g'}} \mc{D}_{f_{g'}}(M)} 
\end{array}
$$
is a constant.
 
Suppose now that $g'$ is a Yamabe metric so
$s_{g'}= n(n-1) + \| H_{f_{g'}}\|^2 - \| \alpha_{f_{g'}}\|^2$ with
$\| H_{f_{g'}}\|^2$ and $\| \alpha_{f_{g'}}\|^2$ constant, and no other 
metric in $[g]$ such that $\mu_{g'}=\mu_g$ has smaller constant 
scalar curvature. 
We consider another Yamabe metric in the conformal class, which 
for convenience we label $g$, of the same volume as $g'$, and assume that 
$\mc{W}_{f_{g}} \geq \mc{W}_{f_{g'}}$, and $\mc{D}_{f_{g}}\geq 
\mc{D}_{f_{g'}}$, respectively. In this setting, we use the path above
connecting $f_{g}$ and $f_{g'}$, and by the equality of the Yamabe functional
values of $g$ and $g'$, we have that
$$
\mc{W}_{f_{g'}}-\mc{D}_{f_{g'}} =\mc{W}_{f_{g}} -\mc{D}_{f_{g}} \, ,
$$
If $s_{g'}=0$, then by integration of (\ref{gra}) against the measure 
$e^{(n-2)u(1)}d\mu_g$, we conclude that $s_{g}=0$ also. If we then
integrate the resulting identity (\ref{gra}) against the measure 
$e^{nu(1)}d\mu_g$, we conclude that $du^{\tau}\mid_{t=1}=0$, and so 
$u(t)\mid_{t=1}$ is constant, which by the equal volume condition, must be 
zero. So $g'=g$ up to isometries,
and the Yamabe metric $g$ is such that
$\mc{W}_{f_{g}}=\mc{W}_{f_{g'}}$, and 
$\mc{D}_{f_{g}}=\mc{D}_{f_{g'}}$. 
If $s_{g'}\neq 0$, 
the quotient of the difference of the critical point equations for $f_g$ and 
$f_{g'}$ yields that
$$
\frac{e^{2u(1)}s_{g'}}{s_g}=e^{2u(1)}=
\frac{\mc{W}_{f_{g'}}-\mc{D}_{f_{g'}}}{\mc{W}_{
f_{g}} -\mc{D}_{f_{g}}}=1 \, ,
$$
and so $u(1)=0$, $g'=g$ up to isometries again, 
and the Yamabe metric $g$ has  
$\mc{W}_{f_{g}}=\mc{W}_{f_{g'}}$, and $\mc{D}_{f_{g}}=\mc{D}_{f_{g'}}$. 
Thus, the only constant scalar curvature metrics that are Yamabe metrics 
in their conformal classes, are the minimizers $f_{g'}$ of the 
functionals $\mc{W}_{f_{g}}(M)/\mu_g^{2/N}$ and $\mc{D}_{f_{g}}(M)/\mu_g^{2/N}$.

A scalar flat or Einstein metric is a Yamabe metric in its class. If $g'$
is any one such, and $f_{g'}$ is not minimal, we consider a path of
conformal homothetics defined by a function $u$ such that
$u(s)\mid_{f_{g'}(M)}=s$ and $(\nabla^{g'} u)\mid_{f_{g'}
(M)}=(\nabla^{g'} u)^{\nu} \mid_{f_{g'}(M)}=H_{f_{g'}}$. Along the 
associated path $[0,1/n] \ni s \rightarrow f_{e^{2s}g'}$,
we have that
$$
\begin{array}{rcl} 
\| H_{f_{e^{2s}g'}}\|^2 & = & e^{-2s}(1-ns)^2\| H_{f_{g'}}\|^2 
\, , \vspace{1mm} \\
n^2+ \| H_{f_{e^{2s}g'}}\|^2 & = & e^{-2s}\left( n^2 -(n(ns^2-2s)
-(1-ns)^2)\| H_{f_{g'}}\|^2 \right)=
 e^{-2s}\left( n^2 +\|H_{f_{g'}} \|^2 \right) \, , \vspace{1mm} \\
n+\| \alpha_{f_{e^{2s}g'}}\|^2 & = & e^{-2s}\left(n-((ns^2-2s)-(2s-ns^2))
\| H_{f_{g'}}\|^2) +\| \alpha_{f_{g'}}\|^2\right) =
e^{-2s}\left(n + \|\alpha_{f_{g'}}\|^2\right)  \, .
\end{array}
$$
At $s=1/n$, $\| H_{f_{e^{2s}g'}}\|^2 = 0$. 
On the other hand, if $g'$ is any Yamabe metric that is neither scalar flat 
nor Einstein, then $f_{g'}$ must be minimal. For otherwise,
since 
$$
\frac{d}{dt} \lambda(M,g'(t))\mid_{t=0} = 
\frac{1}{\mu_{g'}^{\frac{2}{N}}} \int -(r_{g'},h) d\mu_{g'}
+\frac{1}{\mu_{g'}^{\frac{2}{N}}}
s_{g'}\left(1-\frac{2}{N}\right) \frac{d}{dt} \int d\mu_{g'(t)}\mid_{t=0} 
\, , 
$$
where $g'(t)$ is a path of conformal deformations of $g'$ that infinitesimally
varies in the direction of the symmetric two tensor $h=\dot{g}'(0)$, we could
choose a variational path of $f_{g(t)}$, where the metric $\tg$ on the
background sphere $\mb{S}^{\tn}$ is varying in the conformal direction 
$\|H_{f_{g'}}\|\tg(\, \cdot \, , \, \cdot \,)\mid_{f_{g'}(M)}$ at $t=0$, to 
obtain that
$$
\frac{d}{dt} \lambda(M,g'(t))\mid_{t=0} = 
-s_{g'}\left(\frac{1}{\mu_{g'}^{\frac{2}{N}}} \| H_{f_{g'}}\| 
+\frac{1}{\mu_{g'}^{\frac{2}{N}}}
\left(1-\frac{2}{N}\right)\| H_{f_{g'}}\|^2 \right)\, ,  
$$
which is nonzero, and so $g'$ could not be Yamabe. Conversely, if 
$s_{g'}$ is a nonzero constant, and $f_{g'}$ is minimal, 
by \cite[Theorem 4]{sim5}, or by the uniqueness of Yamabe metrics in
the nonpositive case, $g'$ is a Yamabe metric. In all of
these three cases, by the 
homothetic invariance of the functionals $\mc{W}_{f_{g}}(M)/\mu_g^{2/N}$ and 
$\mc{D}_{f_{g}}(M)/\mu_g^{2/N}$, the minimal $f_{c^2g'}$ is a minimizer of
$\mc{W}_{f_{g}}(M)$ and $\mc{D}_{f_{g}}(M)$ in the conformal class.
\qed
\bigskip

\begin{theorem} \label{thM2}
Any smooth path $[0,1]\ni t\rightarrow [g_t]$ 
of conformal classes, which corresponds to a smooth path $t \rightarrow g_t
\in \mc{M}(M)$ of metrics, admits a lift to a 
path $[0,1]\ni t\rightarrow g_t^Y$ of Yamabe metrics $g_t^Y \in [g_t]$
such that $\mu_{g_t^Y}(M)=\mu_{g_t}(M)$, with the paths 
$t\rightarrow \| H_{f_{g_t^Y}}\|^2$ and $t\rightarrow \| 
\alpha_{f_{g_t^Y}}\|^2$ of constant functions being smooth and  
{\rm (}at least{\rm )} continuous,
respectively. If $[g_{t}]=[g_n]$ for at least one $t \in [0,1]$,
then $M^n$ is diffeomorphic to $\mb{S}^n$, 
$[g_t]=[g_n]$ for all $t\in [0,1]$, and 
$g^Y_t=\left( \frac{\mu_{g_t}}{\omega_n} \right)^{\frac{2}{n}}\tg\mid_{f_{g_n}
(\mb{S}^{n})}$.   
\end{theorem}

{\it Proof}. In the $C^2$ topology on the space of metrics the 
mapping $g\rightarrow \lambda(M,[g])$ is continuous 
\cite[Proposition 7.2]{bebe}, hence, in the resulting quotient topology in 
the space of conformal classes, the mapping $[g]\rightarrow \lambda(M,[g])$ is 
continuous also. We can thus partition the interval $[0,1]$ into 
$\{t_0=0, \ldots, t_p=1\}$ if necessary so that, on every subinterval 
$[t_j,t_{j+1}]$,  
either $\lambda(M,[g_t])\leq 0$, or
$\lambda(M,[g_t]) \geq 0$ and $>0$ on $(t_j,t_{j+1})$.
By the uniqueness of the solution to the Yamabe problem on any conformal class 
of metrics on $M$ in the nonpositive case, the result follows on each of 
the subintervals where $\lambda(M,[g_t])\leq 0$, and therefore, by
the trychotomy theorem of Kazdan \& Warner \cite{kawa}, the result follows
for any manifold $M$ for which not every smooth function on it is the scalar 
curvature of some smooth metric.  In these cases, the resulting path 
$t\rightarrow s_{g_t^Y}$ of scalar curvatures of the Yamabe metrics 
$g_t^Y$ is differentiable. This reduces the proof of the Theorem to the 
case where $M$ is assumed to be a manifold that carries metrics of positive 
scalar curvature, and under that assumption, for curves $t\rightarrow [g_t]$ 
such that $[0,1]\ni t\rightarrow \lambda(M,[g_t])\geq 0$, and $>0$
for at least one $t\in (0,1)$. 

By Theorem \ref{thM1}, a Yamabe metric $g_t^Y$ in $\mc{M}_{[g_t]}(M)$ is 
either scalar flat, or Einstein of positive scalar curvature, or otherwise it
is of constant positive scalar curvature and has
minimal isometric embedding $f_{g_t^Y}$. By the uniqueness of the solution
to the Yamabe problem in the nonpositive case, in the first case, the 
scalar flat Yamabe metric in $[g_t]$ of volume $\mu_{g_t}$ is unique; 
by Obata's theorem \cite{oba}, in the second case, the Einstein Yamabe metric 
in $[g_t]$ of volume $\mu_{g_t}$ is unique; by the second identity in 
(\ref{pti}), two equal volume Yamabe metrics in the same class with minimal 
isometric embeddings are equal to each other, so in the third case, the Yamabe 
metric with minimal isometric embedding, and volume $\mu_{g_t}$,  
is unique also. Since the path $t\rightarrow \lambda(M,[g_t])$ is continuous,
by lifting to these unique Yamabe metrics in the classes, 
we produce a path $t \rightarrow f_{g_t^Y}\in 
\mc{M}_{\mb{S}^{\tn}}(M)$ of isometric embeddings of Yamabe metrics with the 
desired properties. Along this path, $t\rightarrow s_{g_t^Y}$ is  
(at least) continuous, and the three nonintrinsic 
functionals $\Theta_{f_g}(M)$, 
$\Phi_{f_g}(M)$ and $\Pi_{f_g}(M)$ have densities given by the paths of 
constant functions  $t\rightarrow n(n-1)$, $t\rightarrow \| H_{f_{g_t^Y}}\|^2$ 
and $t\rightarrow \| \alpha_{f_{g_t^Y}}\|^2$, respectively. The first and 
second of these paths vary differentiably, and so the third must vary
at least continuously. Indeed, by the Palais isotopic extension theorem, 
we may produce a smooth deformation $f_{g_{s,t}}$ of $f_{g_t}$ parametrized by 
a smooth two parameter family of functions $[0,1]\times [0,1] \ni (t,s) 
\rightarrow u_t(s)$ such that $g_{s,t}=e^{2u_{t}(s)} \tg$ on $f_{g_{s,t}}(M)$, 
$g_{0,t}=g_t$, $[g_{s,t}]=[g_t]$, $\mu_{g_{s,t}}=\mu_{g_t}$, and 
$g_{1,t}=g_t^Y$, so that if $f_{g_{s,t}}$ is the associated two parameter 
family of isometric embeddings of $g_{s,t}$, then 
$t\rightarrow f_{g_{1,t}}=f_{g_t^Y}$. 
By the second of the identities in (\ref{pti}), $\| H_{f_{g_{s,t}}}\|^2$ 
depends differentiably on $\| H_{f_{g_t}}\|^2$, $u_t(s)$, and the normal 
component of $\nabla^{\tg} u_t(s)$ (in fact, by 
\cite[Theorem 3.10]{gracie}, we have that this constant function
satisfies the elliptic equation
$$
\Delta^{g_{s,t}} \left( \| H_{f_{g_{s,t}}}\|^2 \right) = 0 =
\| H_{f_{g_{s,t}}}\|^2 \left( 2n- 2\| \nabla_{e_i}^\nu \nu_{
H_{f_{g_{s,t}}}}\|^2 - \| H_{f_{g_{s,t}}}\|^2 + 2 {\rm trace}\,A_{\nu_{
H_{f_{g_{s,t}}}}}^2\right)\, , 
$$
where $\nu_{H_{f_{g_{{s,t}}}}}$ is a normal vector in the direction of
$H_{f_{g_{s,t}}}$, $A_{\nu_{H_{f_{g_{{s,t}}}}}}$ is the shape operator along 
that direction, and $\nabla^{\nu}$ is the covariant derivative of the normal 
bundle of the embedded submanifold $f_{g_{s,t}}(M)$).  
Hence, $\| H_{g_{1,t}}\|^2=\| H_{f_{g_t^Y}}\|^2$ depends differentiably on $t$.

Suppose now that $[g_t]=[g_n]$ for some $t=t_0 \in [0,1]$, and so
we have that $\lambda(M,[g_{t_0}])=\lambda(\mb{S}^{n},[g_n])$, and 
that $[g_{t_0}]=[g_{t_0}^Y] =[g_n]$.
Hence, $f_{g_{t_0}}(M)$ is conformally isotopic to the standard geodesic 
embedding $f_{g_n}(\mb{S}^n,g_n)\hookrightarrow (\mb{S}^{\tn},\tg)$, which
implies that $M \cong \mb{S}^{n}$, and by the solution of the Yamabe problem 
on the standard sphere, $g^Y_{t_0}$ is, up to a dilation factor, obtained 
from $g_n$ by a conformal diffeomorphism of $\mb{S}^n$ \cite{oba}. Thus,  
$g^Y_{t_0}=e^{2u_{t_0}(s)\mid_{s=1}}\tg\mid_{f_{g_n}(\mb{S}^n)}$ for a certain 
path of functions $[0,1] \ni s\rightarrow u_{t_0}(s)$ such that $u_{t_0}(0)=0$, 
relating conformally to each other the isometric embeddings
$f_{g_n}$ and $f_{g^Y_{t_0}}(M)$. Since 
either $g^{Y}_{t_0}$ is Einstein, or $f_{g^Y_{t_0}}$ is minimal, 
and since by \cite[Theorem 5.2.3]{simo}, there 
exists a $C^2$-neighborhood $\mc{E}'$ of $f_{g_n}$ in the space of
of smooth immersions of $\mb{S}^{n}$ into $\mb{S}^{\tn}$ such that no 
immersion in it other than $f_{g_n}$ itself is minimal, we conclude that 
$g^Y_{t_0}$ must be a homothetic deformation of $g_n$, hence Einstein, 
with dilation factor determined by $\mu_{g_{t_0}}$. We then   
consider an $\varepsilon$ neighborhood 
$U^{\varepsilon}_{t_0}$ of $t_0$ in $[0,1]$,
and the restriction of the path of Yamabe metrics to this neighborhood. By the 
Palais isotopic extension theorem, we produce a family of 
functions $U^{\varepsilon}_{t_0} \times [0,1] \ni (t,s) \rightarrow u_{t}(s)$ 
such that $U^{\varepsilon}_{t_0} \ni t \rightarrow 
f_{g_t^{Y}}=f_{e^{2u_t(s)\mid_{s=1}}\tg\mid_{f_{g_n}(\mb{S}^n)}}$, 
and with the family of isometric embeddings
$U^{\varepsilon}_{t_0} \ni t \rightarrow 
f_{e^{2u_t(s)\mid{s=0}}\tg\mid_{f_{g_n}(\mb{S}^n)}}$ all contained in some  
neighborhood $\mc{E}$ of $f_{g_n}$ in $\mc{M}_{\mb{S}^{\tn}}(M)$. We choose
$\varepsilon=\varepsilon(g_n)>0$ small enough such that 
$\mc{E}\subset \mc{E}' \cap \mc{M}_{\mb{S}^{\tn}}(M)$. Since $g_t^Y$ is
either Einstein, or $f_{g_t^Y}$ is minimal, 
it follows that 
$U^{\varepsilon}_{t_0} \ni t\rightarrow g_t^Y = 
e^{2u_t(s)\mid_{s=1}}\tg_{f_{g_n}(\mb{S}^n)}$ 
must be a path of dilations of $g_n$, and so all the $u_t(s)\mid_{s=1}$ 
are constant functions, and $[g_t]=[g_t^Y]=[g_n]$.
Since the choice of $\varepsilon$ depends on $g_n$
only, and the functions $[0,1]\ni t \rightarrow [g_t]$ and
$g\rightarrow \lambda(M,[g])$ are continuous, we can iterate this argument
to prove that $g_t^Y$ is a smooth dilation of $g_n$, and that 
$[g_t]=[g_t^Y]=[g_n]$ for all $t \in [0,1]$. 
The remaining part of the proof concerning the dilation factors as 
a function of $\mu_{g_t}$ and $\omega_n$ is 
straightforward.
\qed

\subsection{Metrics of constant $J$ scalar curvature: $J$
Yamabe metrics}
If a manifold $M^{n=2m}$ carries an almost complex structure $J_0$, 
we denote by $\mc{J}_{J_0}(M)$ the manifold of all almost complex structures 
on $M$ in the same orientation class as $J_0$. If 
$J\in \mc{J}_{J_0}(M)$, a metric $g$ is Hermitian relative to 
$J$ if $g(JX,JY)=g(X,Y)$, in which case we say that $J$ is compatible  with $g$,
and that $g$ is invariant under $J$.
We denote by $\mc{I}_{J_0}(M,g)$ the space of all almost complex 
structures in $\mc{J}_{J_0}(M)$ compatible with $g$.
We define    
$\mc{M}^{J_0}(M)$ to be the space of metrics on $M$ that are compatible with
at least one element of $\mc{J}_{J_0}(M)$, and $\mc{C}^{J_0}(M)$ to be the
space of conformal classes of such metrics. We then obtain a fibration
\begin{equation} \label{eq20}
\mc{M}^{J_0}(M) \stackrel{\pi}{\rightarrow } \mc{C}^{J_0}(M) \, .  
\end{equation}
The set $\mc{M}^{J_0}_{[g]}(M)$ of metrics in the conformal class $[g]$ 
of $g \in \mc{M}^{J_0}(M)$ is its fiber over $[g]$, 
and we have the foliated decomposition $\mc{M}^{J_0}(M)=\cup_{[g]\in
 \mc{C}^{J_0}(M)}\mc{M}^{J_0}_{[g]}(M)$. 
Any metric $g$ in $[g]\in \mc{C}^{J_0}(M)$ is invariant under any element of 
$\mc{I}_{J_0}(M,g)$.

If $g\in \mc{M}^{J_0}(M)$, we choose a $J\in \mc{I}_{J_0}(M,g)$, and define the
$J$-scalar curvature $s^J_{g}$ of $g$ as the metric trace of the $J$-Ricci 
tensor 
$$
r_{g}^J(X,Y)={\rm trace}\, L \rightarrow -J(R^{g}(L,X)JY) \, .
$$
This $J$-scalar curvature function is related to the scalar curvature of $g$ 
by the identity
\begin{equation} \label{rem} 
(n-1)s_g^J-s_g=2(n-1)W_g(\omega_g^{\sharp},\omega_g^{\sharp})\, , 
\end{equation}
where $W_g$ is the Weyl tensor of $g$, and $\omega_g^{\sharp}$ is the 
fundamental form viewed as a bivector \cite[Proposition 6]{rss}, and 
though the skew Hermitian component of $r^J_g$ is a conformal 
invariant that may very with $J$ \cite[Corollary 4.5]{rss},  
$s_g^J$ depends on $J\in \mc{I}_{J_0}(M,g)$ only
through the rotational component of the polar decomposition of the triple
$(g,J,\omega_g)$ on the tangent space at each point, 
and so it is a function of $\mc{I}_{J_0}(M,g)$ but not of the particular 
choice of $J$ in it made to define it. Thus, the $J$-Yamabe 
functional  
\begin{equation} \label{jyfu}
\mc{M}^{J_0}_{[g]}(M) \ni g \rightarrow \lambda^J(g) = 
\frac{1}{\mu_g^{\frac{2}{N}} }\int s_g^J d\mu_{g} 
\end{equation}
is well-defined.

\begin{theorem} \label{thM3}
If $g_t = e^{2u(t)}g$ is a path of metrics in the conformal class of the $J$
invariant metric $g$, of associated isometric embeddings $f_{g_t}$, 
$g'=g_t{\mid_{t=1}}$ is a critical point of 
{\rm (\ref{jyfu})} if, and only if, $s^J_{g'}$ is a constant satisfying the
equation  
$$
s^J_{g'}= e^{-2u}\left( \frac{1}{n-1}s_g +2
W_{g}(\omega^{\sharp}_g, \omega^{\sharp}_g)  
+ 2 \Delta^{g} u^{\tau} -(n-2)\tg(du^{\tau},du^{\tau}) 
\right)_{\mid_{t=1}} 
= \frac{1}{\mu_{g'}}\int_{f_{g'}(M)} \hspace{-0.1cm} s^J_{g'} d\mu_{g'}\, . 
$$
Such a metric has constant scalar curvature $s_{g'}$ also if, and only if, 
$W_{g'}(\omega^{\sharp}_{g'}, \omega^{\sharp}_{g'})$ is a constant function,
in which case, $(n-1)\lambda^J(g')=\lambda(g')+2(n-1)W_{g'}( \omega^{\sharp}_{
g'}, \omega^{\sharp}_{g'})\mu_{g'}^{\frac{2}{n}}$.
\end{theorem}

{\it Proof}. The variation of (\ref{jyfu}) in the space of deformations of 
the isometric embeddings of the metrics $g_t=e^{2u}g$ is
given by 
$$ 
\begin{array}{rcl}
{\displaystyle \frac{d}{dt} \frac{1}{\mu^{\frac{2}{N}}_{g_t}}
\int_{f_{g_t}(M)}s_{g_t}^J d\mu_{g_t} } & = & 
{\displaystyle  \frac{n-2}{\mu^{\frac{2}{N}}_{g_t}} \left(
\int_{f_g(M)} e^{nu}\dot{u}\left( e^{-2u}\left( \frac{1}{n-1}s_g +2
W_{g}(\omega^{\sharp}_g, \omega^{\sharp}_g)  
\right. \right. \right. } \vspace{1mm} \\ & & 
- {\displaystyle \left. \left. \left. 
(n-2)\tg(du^{\tau}, du^{\tau})
+ 2 \Delta^{g} u \right) 
-\frac{1}{\mu_{g_t}}\int_{f_{g_t}(M)} s^J_{g_t}
  \right) d\mu_{g_t}  \right) }    
\, .  
\end{array}
$$
The desired result follows. 
\qed

Metrics that realize the intrinsic conformal invariant 
$$ 
\lambda^{J}(M,[g])= \inf_{g\in \mc{M}^{J_0}_{[g]}(M)} 
\frac{1}{\mu_g(M)^{\frac{2}{N}}}\int s^J_g d\mu_g  
$$ 
exists \cite{rss}, and are by definition the $J$ Yamabe metric representatives 
of $[g]$. By Aubin's \cite[Theorem p. 155]{au2}, we have the universal bound
$$
\lambda^{J}(M,[g])\leq 2m \omega_{2m}^{\frac{1}{m}} \, 
= n \omega_{n}^{\frac{2}{n}} \, ,
$$
and since the left side varies continuously as a function of $J,[g]$ \cite[
Lemma 2]{sim5}, 
\begin{equation} \label{sigJ0}
\sigma^{J_0}(M)= \sup_{[g]\in \mc{C}^{J_0}(M)} \lambda^J(M,[g]) 
\end{equation}
is a well-defined invariant of the differentiable almost Hermitian 
structure on $M$. 

The following is now straightforward.

\begin{theorem} \label{th5}
{\rm \cite[Theorem 3]{sim5}} A metric $g$ in $\mc{M}^{J}(M)$ 
solves both, the Yamabe and $J$ Yamabe problem in its class, if, 
and only if, $W_g(\omega^{\sharp}_g, \omega^{\sharp}_g)=0$. 
\end{theorem}

By identity (\ref{rem}), the solutions to the Yamabe and almost Hermitian 
Yamabe problems lead to the following result. If combined with Simons' 
gap theorem, it serves to single out the octonionic triple $(\mb{S}^6,J,g)$, a
fact tacitly present already in the proof of Theorem \ref{th1} (a)(c).

\begin{theorem}
Any smooth path $[0,1]\ni t\rightarrow [g_t]\in \mc{C}^{J_0}(M)$, which
corresponds to a smooth path $t \rightarrow (J_t,g_t) \in 
\mc{I}_{J_0}(M)\times \mc{M}^{J_0}(M)$ of $J_t$ invariant metrics $g_t$, 
admits a lift to a path $[0,1]\ni t\rightarrow g_t^{JY}$ of $J$ Yamabe metrics 
of the form $g_t^{JY}=e^{2u_t}g_t^Y$, $g_t^Y$ a Yamabe metric in $[g_t]$,
and $\mu_{g_t^{JY}}(M)=\mu_{g_t^Y}(M)=\mu_{g_t}(M)$.
If $[g_{t}]=[g_n]$ for at least one $t \in [0,1]$, then $M^{n=6}$ is 
diffeomorphic to $\mb{S}^6$, $J_t \in\mc{I}_{J}(\mb{S}^6,g)$, 
$[g_t]=[g_6]$ for all $t\in [0,1]$, and  
$f_{g_t^{JY}}\in \mc{M}_{\mb{S}^{\tn}}(\mb{S}^6)$ is a homothetic deformation 
of the totally geodesic embedding $f_{g}$ of $\mb{S}^6$ corresponding to the 
metric $g_t^{JY}=(\mu_{g_t}/\omega_6)^{\frac{1}{3}}
\tg_{\mid_{f_{g}(\mb{S}^6)}}$.
\end{theorem}

\subsection{Product of spheres with $\sqrt{-\BOne}$ tensors on their tangents}
If $n,k\in \mb{N}$, $k<n$, we consider the manifold 
$M^n_k=\mb{S}^k \times \mb{S}^{n-k}$, and if $r\in (0,1)$, we set
$$
M^{n}_k(r,\sqrt{1-r^2})= \mb{S}^k(r) \times \mb{S}^{n-k}(\sqrt{1-r^2})\, ,
$$
and let $g^n_k=g^n_k(r)$ be the product metric on it. We view  
$(M^n_k(r,\sqrt{1-r^2}),g^n_k)$ isometrically embedded into 
$(\mb{S}^{n+1},g) \hookrightarrow (\mb{S}^{\tn},\tg)$. 
We distinguish the special case
\begin{equation} \label{mipr} 
\mb{S}^{n,k} = M^n_k\left(\sqrt{\frac{k}{n}},\sqrt{\frac{n-k}{n}}\right)
\subset \mb{S}^{n+1} \, ,
\end{equation}
and denote its product metric by $g_{\mb{S}^{n,k}}$. 
The isometric embedding $f_{g_{\mb{S}^{n,k}}}$ is minimal, and
$s_{g_{\mb{S}^{n,k}}}=n(n-2)$, so the
metric $g_{\mb{S}^{n,k}}$ is a Yamabe metric in its conformal class
\cite[Theorem 4, Corollary 5]{sim5} (cf. Theorem \ref{thM1}). 

If $n=2(l+m+1)$ is even and $k=2l+1$ is odd, $M^n_k(r,\sqrt{1-r^2})$ is
the product of two odd dimensional spheres, an example of a Calabi-Eckmann 
manifold carrying an integrable complex structure $J^r_{n,k}$ compatible with 
$g^n_k$ \cite{caec}. When $r=\sqrt{(2l+1)/n}$, we refer to the complex structure
simply by $J_{n,k}$. The metric $g_{\mb{S}^{n,k}}$ has $J_{n,k}$ scalar 
curvature $s^{J_{n,k}}_{g_{\mb{S}^{n,2l+1}}}=(2l/k+2m/(n-k))n$, and is
a $J_{n,k}$ Yamabe metric in its conformal class if, only
if, $ml=0$ \cite[Example 10 (2)]{sim5} (if $n=2$, $l=m=0$,  
$M^n_k(\sqrt{1/2},\sqrt{1/2})$ is the Clifford torus, and  
$s^{J_{n,k}}_{g_{\mb{S}^{n,k}}}=s_{g_{\mb{S}^{n,k}}}=0$, so $g^n_k$ is Ricci
flat).

The product $\mb{S}^{4,2}$ is naturally K\"ahler with an integrable 
almost complex structure $J$. There are octonionic induced almost complex 
structures $J$ on each of ${\mb{S}^{6,2}}$, ${\mb{S}^{8,2}}$, and 
${\mb{S}^{12,6}}$, respectively, compatible with $g_{\mb{S}^{n,k}}$, the 
latter two of these splitting the manifold as a product of almost Hermitian 
spheres. In any of these cases, the metric $g_{\mb{S}^{n,k}}$ is not a $J$ 
Yamabe metric in its conformal class \cite[Example 10 (3)]{sim5}.

\begin{theorem}
{\rm \cite{caec} (a)} Products of odd dimensional spheres admit integrable 
almost complex structures. {\rm \cite{dasu} (b)} No product of even dimensional
spheres of dimensions $u$ and $v$, with $\min{\{ u,v\}}>0$, 
admits almost complex structures unless the set of dimensions is any of 
$\{2,2\}$, $\{2,4\}$, $\{2,6\}$, or $\{6,6\}$, respectively.
{\rm \cite[Theorem 12]{sim5} (c)} None of the almost complex structures on 
any of the manifolds $\mb{S}^2 \times \mb{S}^4$, $\mb{S}^2 \times \mb{S}^6$, 
or $\mb{S}^6 \times \mb{S}^6$, respectively, is integrable.
\end{theorem}

{\it Proof}. (a) If the product involves at least one exotic factor, we 
may use an exotic diffeomorphism of this and a Calabi-Eckmann product to 
identify tensors in one with tensors in the other. The integrable 
Calabi-Eckmann almost complex structure on the standard product gets 
identified with an integrable almost complex structure on the exotic one.

(b) Suppose that we have a product of standard spheres, and let $J'$ 
be an almost complex structure on the manifold 
$\mb{S}^{n,k}$, where $n=u+v$, and $k=\min{\{u,v\}}$. If 
$g=g_{\mb{S}^{n,k}}$,
we consider the smooth path of metrics 
$$ 
[0,1]\ni t \rightarrow g_t( \, \cdot \, , \, \cdot \, ) = 
(1-t) g(J' \, \cdot \, , J' \, \cdot \, ) + 
t g( \, \cdot \, , \, \cdot \, ) \, ,  
$$ 
and obtain a path of metrics on $M$ that begins at $g_0=g\circ J'$, and ends 
at $g$. We let $f_{g_t}: 
(M,g_t) \rightarrow (\mb{S}^{\tn},\tg)$ be the associated path of Nash 
isometric embeddings in $\mc{M}_{\mb{S}^{\tn}}(M)$. 
We have that $g_1=g$, and so $[g_1]=[g]$. We show that $[g_t]=[g]$ for all $t$. 

By Theorem \ref{thM2}, there is a unique path of isometric 
embeddings of Yamabe metrics $g_t^Y$ in $[g_t]$ such that $\mu_{g_t^Y}=
\mu_{g_t}$, and ${g^Y_t \mid_{t=1}}=f_g$. We normalize and
rename these metrics to make them all of volume $\mu_g(M)$. If 
$$
f_{g^{Y}_t}: (M,g_t^Y) \hookrightarrow  (\mb{S}^{\tn},\tilde{g}) 
$$
is the associated path of isometric embeddings of the path of Yamabe 
metrics, we have that
$$
s_{g_t^Y}=\frac{1}{\mu_g^{\frac{2}{n}}}\lambda(M,[g_t])=
n(n-1) + \| H_{f_{g^{Y}_t}}\|^2 -\| \alpha_{f_{g^{Y}_t}} \|^2 \, , 
$$
and the paths of constant functions $t\rightarrow s_{g_t^Y}$ and $t\rightarrow 
\| \alpha_{f_{g^{Y}_t}}\|^2$ are (at least) continuous, and 
$t\rightarrow \| H^Y_{f_{g^{Y}_t}}\|^2$ is smooth. 
We then may dilate the metric $g_t^Y$ on $f_{g_t^Y}(M)$ in the direction of 
$H_{f_{g_t^Y}}$, to obtain a family of Yamabe metrics 
$\tilde{g}_t^Y= (1+\| H_{f_{g_t^Y}} \|^2/n^2)g_t^Y$ in $[g_t]$, whose family of
isometric embeddings
$$
f_{\tilde{g}_t^Y}: (M,\tg_t^Y) \rightarrow (\mb{S}^{\tn},\tilde{g})
$$
is minimal for all $t$, and the paths
of constant functions $t\rightarrow s_{\tg_t^Y}$ and $t\rightarrow 
\| \alpha_{f_{\tg^{Y}_t}}\|^2$ are (at least) continuous.
We have $\tilde{g}_t^Y \mid_{t=1}=g_t^Y\mid_{t=1}=g$,
and $f_{\tilde{g}_1^Y}=f_{g_1^Y}=f_{g}$.

At $t=1$, we have that $\| \alpha_{f_{\tilde{g}_1^Y}}\|^2=n$,  
$\| H_{f_{\tilde{g}_1^Y}}\|^2=0$, $f_{\tilde{g}_1^Y}=
f_{g}$, and $[\tg_1^Y]=[g_1^Y]=[g]$. By the reinterpretation of the gap 
theorem of Simons \cite[Corollary 5.3.2]{simo} in terms of critical points of 
$\Psi_{f_g}(M)$ in the space of isometric embeddings of metrics in the 
conformal class \cite[Theorem 2]{rss2}, we conclude that 
modulo isometries, this remains true for any 
$t$ for which $\| \alpha_{f_{\tilde{g}_t^Y}}\|^2$ is sufficiently close
to $n$ \cite[Theorem 9]{sim5}. Thus, there exists a smallest 
$a$, $0\leq  a \leq 1$ such that if $t\in [a,1]$, $f_{\tilde{g}_t^Y} = 
f_{g_t^Y}=f_{g}$, $\| \alpha_{f_{\tg_t^Y}}\|^2=n$, while for $t<a$ we would
have that $\| \alpha_{f_{\tilde{g}^Y_t}}\|^2 > n$. By the 
continuity of the entire path $t \rightarrow \| \alpha_{f_{\tg^{Y}_t}}\|^2$, 
the alluded theorem implies that this smallest $a$ is $a=0$. 
Thus, $[\tg_t^Y]=[g_t]=[g]$ for all $t$.  

Since the path of classes $t \rightarrow [g_t]$ is constant, $[g_0]=[g]$, and 
$g_0=g\circ J'$ and $g$ are conformally related. This implies that $g$
is $J'$-invariant, and that the pair $(J',g)$ is an almost Hermitian 
structure on $\mb{S}^{n,k}$. If we consider the linear isometric embedding
of $(\mb{S}^{n,k},g_{\mb{S}^{n,k}})$ into $(\mb{S}^{n+1},g)$, if the
almost Hermitian pair $(J',g)$ does not split, by Kirchhoff theorem 
\cite{kir}, $\mb{S}^{n+1}$ must be parallelizable, so
$n=6$, and the manifold must be $\mb{S}^{6,2}$ with an octonionic 
almost complex structure $J$.  Otherwise, the factors of the manifold are 
spheres with almost complex structures, and by Theorem \ref{th1}, it 
must then be either $\mb{S}^{4,2}$ with its product K\"ahler structure $J$, or 
$\mb{S}^{8,2}$ or $\mb{S}^{12,6}$ with their octonionic induced $J$s. In all of 
these cases, we have that
$J'\in \mc{I}_{J}(\mb{S}^{n,k},g_{\mb{S}^{n,k}})$.

If we start now with a product of spheres where at least one of the factors
has an exotic differentiable structure, we use an associated exotic 
diffeomorphism between this product and $\mb{S}^{n,k}$, and by pull-back and 
push-forward, identify almost complex structures in the former with almost 
complex structures in the latter. The proof above then carries onto the 
product of spheres with at least one exotic factor.

(c) We let $(M^n,J,g)$ stand here for any of the almost Hermitian manifolds 
$(\mb{S}^{6,2},J,g_{\mb{S}^{6,2}})$, $(\mb{S}^{8,2},J,g_{\mb{S}^{8,2}})$, and 
$(\mb{S}^{12,6},J,g_{\mb{S}^{12,6}})$ with their octonionic almost complex
structures. We have that $g$ is a Yamabe metric in its conformal class,
$s_g=n(n-2)$, and $W_{g}(\omega_g^{\sharp},\omega_g^{\sharp})=n^2$. It follows
that 
$$
\lambda(g)=\lambda(M,[g])< \lambda^J(g)=n(n-1)\mu_{g}^{\frac{2}{n}} < 
n(n-1)\omega_n^{\frac{2}{n}}\, ,
$$
and by Theorem \ref{th5}, $g$ is not a $J$ Yamabe metric in its conformal 
class, and $\lambda^J(M,[g])< n\omega_n^{\frac{2}{n}}$.
We choose and fix a $J$ Yamabe metric $g^{JY}$ in the class of volume 
$\omega_n$, and have that 
\begin{equation} \label{ub}
n\frac{n-2}{n-1}\left( \frac{\mu_g(M)}{\omega_n}\right)^{\frac{2}{n}}
< s_{g^{JY}} < n \, .
\end{equation}
The quantity on the left here lies on the interval $(0.91n, 0.97n)$.
 
We consider any $J'\in \mc{I}_{J}(M)$, and choose a smooth path
$[0,1]\ni t \rightarrow J_t$ of almost complex structures connecting 
$J'$ and $J$. We define the path of 
smooth metrics 
$$
[0,1]\ni t \rightarrow g_t (\, \cdot \, , \, \cdot \,)=\frac{1}{2}(
g(\, \cdot \, , \, \cdot \, ) + g(J_t \, \cdot \, , J_t \, \cdot \, ))\, ,
$$
and produce a path $t \rightarrow (J_t,g_t)$ of almost Hermitian structures on
$M$ that begins at $(J',g_0)$, and ends at $(J,g)$. We let $f_{g_t}: 
(M,g_t) \rightarrow (\mb{S}^{\tn},\tg)$ be the associated path of Nash 
isometric embeddings in $\mc{M}_{\mb{S}^{\tn}}(M)$.

By Theorem \ref{thM2}, there is a unique path of isometric 
embeddings of Yamabe metrics $g_t^Y$ in $[g_t]$ such that $\mu_{g_t^Y}=
\mu_{g_t}$, and ${g^Y_t \mid_{t=1}}=f_g$. We normalize and
rename these metrics to make them all of volume $\mu_g(M)$. We also
solve the $J_t$ Yamabe problem in each almost Hermitian
class $(J_t, [g_t])$, and produce a path $g^{JY}_t$ of $J_t$ Yamabe metrics, 
normalized each to have volume $\omega_n$, which at $t=1$ is $g^{JY}$.
We obtain paths of metrics $t \rightarrow g_t^Y$ and 
$t\rightarrow g_t^{J_tY}$ of volume $\mu_g$ and $\omega_n$, 
which at $t=1$ are $g$ and $g^{JY}$, respectively. If
$$
f_{g^{Y}_t}: (M,g_t^Y) \hookrightarrow  (\mb{S}^{\tn= n +p },\tilde{g}) \, , 
$$
and
$$
f_{g^{JY}_t} : (M,J_t,g_t^{JY}) \hookrightarrow  (\mb{S}^{\tn= n +p},
\tilde{g})\, ,
$$
are the associated families of Nash isometric embeddings, we have that
$$
s_{g_t^Y}=\frac{1}{\mu_g^{\frac{2}{n}}}\lambda(M,[g_t])=
n(n-1) + \| H_{f_{g^{Y}_t}}\|^2 -\| \alpha_{f_{g^{Y}_t}} \|^2 \, , 
$$
and 
$$
s_{g_t^{JY}}^{J_t}= \frac{1}{\omega_n^{\frac{2}{n}}} \lambda^{J_t}(M,[g_t]) 
=n+\frac{1}{n-1}(\| H_{f_{g_t^{J_tY}}}\|^2 -\| \alpha_{f_{g_t^{J_tY}}}\|^2)+ 
2W_{g_t^{J_tY}}(
\omega_{g_t^{J_tY}}^{\sharp},\omega_{g_t^{J_tY}}^{\sharp})    
 \, , 
$$
respectively.
The paths of constant functions $t\rightarrow s_{g_t^Y}$ and $t\rightarrow 
\| \alpha_{f_{g^{Y}_t}}\|^2$ are (at least) continuous, while the
path $t\rightarrow \| H^Y_{f_{g^{Y}_t}}\|^2$ is smooth. 

By dilating the metric $g_t^Y$ on $f_{g_t^Y}(M)$ in the direction of 
$H_{f_{g_t^Y}}$, we obtain a family of Yamabe metrics 
$\tilde{g}_t^Y= (1+\| H_{f_{g_t^Y}} \|^2/n^2)g_t^Y$ in $[g_t]$ whose family of
isometric embeddings
$$
f_{\tilde{g}_t^Y}: (M,\tg_t^Y) 
\rightarrow (\mb{S}^{n+p},\tilde{g})
$$
is minimal for all $t$. 
We have $\tilde{g}_t^Y \mid_{t=1}=g_t^Y\mid_{t=1}=g$,
and $f_{\tilde{g}_1^Y}=f_{g_1^Y}=f_{g}$.

At $t=1$, we have that $\| \alpha_{f_{\tilde{g}_1^Y}}\|^2=n$,  
$\| H_{f_{\tilde{g}_1^Y}}\|^2=0$, 
and modulo an isometry, $f_{\tilde{g}_1^Y}=
f_{g_1^Y}=f_{g}$. By Simons' gap theorem \cite[Corollary 5.3.2]{simo} in 
its reinterpreted version \cite[Theorem 9]{sim5}, 
this remains true for any $t$ for which $\| \alpha_{f_{\tilde{g}_t^Y}}\|^2$ is 
nearby $n$. Hence, by the continuity of the path
$t\rightarrow \| \alpha_{f_{\tilde{g}_t^Y}}\|^2$, there exists
a smallest $a$, $0\leq  a \leq 1$ such that if $t\in [a,1]$, 
$f_{\tilde{g}_t^Y}=f_{g}$, and $[\tg_t^Y]=[g]$. 
  
There are no orthogonal complex structures on 
$(M^n,g)$ \cite{euse}. Hence, if we assume that the path $J_t$
starts at an integrable $J'=J_0$, we have that
$[{\tilde{g}^Y_0}]=[g^Y_0]=[g_0]\neq [g]$, and by Simons gap theorem, we must
have $\| \alpha_{f_{\tilde{g}^Y_0}}\|^2>n$. Thus, the $a$ above must be
strictly positive, and so there exists $\varepsilon$, $0< \varepsilon < a$, 
such that for any $t \in (a-\varepsilon, a)$, the functions
$s_{\tilde{g}_t^Y}$, $W_{\tilde{g}_t^Y}(\omega_{\tilde{g}_t^Y}^\sharp,
 \omega_{\tilde{g}_t^Y}^\sharp)$, and $s_{g_t^{J_tY}}^{J_t}$ are positive, 
$[g_t] \neq [g]$, and in this range of $t$s, 
$$
\| \alpha_{f_{\tilde{g}^Y_t}}\|^2 > n \geq \frac{np}{2p-1} \, . 
$$
This leads to a contradiction.

Indeed, since $g_t^{J_tY}$ and $\tg_{t}^{Y}$ are in the same conformal class,
by (\ref{rem}) and the identities in (\ref{pti}), 
there exists a path $t\rightarrow u_t$ of functions such that
$g_t^{J_tY}=e^{2u_t}\tg_t^{Y}$, and 
$$
(n-1)s_{g_t^{JY}}^{J_t}=e^{-2u_t}(n(n-1)-\| \alpha_{f_{\tilde{g}_t^Y}} \|^2+
2(n-1)W_{\tilde{g}_t^Y}(\omega_{\tilde{g}_t^Y}^\sharp, 
\omega_{\tilde{g}_t^Y}^\sharp) + (n-1)(2\Delta^{\tilde{g}_t^Y}u_t - 
(n-2)\tilde{g}_t^Y(\nabla^{\tilde{g}_t^Y}u_t^\tau,\nabla^{\tilde{g}_t^Y}
 u_t^\tau))) \, .
$$
By integration with respect to the measure $e^{\frac{n+2}{2}u_t}d\mu_{
g_t^Y}$,
we obtain that 
$$
(n-1)s_{g_t^{JY}}^{J_t}\int e^{\frac{n+2}{2}u_t}d\mu_{g_t^Y} =
\int (n(n-1)-\| \alpha_{f_{\tilde{g}_t^Y}} \|^2+
2(n-1)W_{\tilde{g}_t^Y}(\omega_{\tilde{g}_t^Y}^\sharp, 
\omega_{\tilde{g}_t^Y}^\sharp))e^{\frac{n-2}{2}u_t}d\mu_{g_t^Y}\, ,  
$$
and since at $t=a$, $\tilde{g}_t^{Y}=g$, $J_t$ and $J$ are compatible with
$g$, and 
$(n(n-1)-\| \alpha_{f_{\tilde{g}_t^Y}} \|^2+ 2(n-1)W_{\tilde{g}_t^Y}(
\omega_{\tilde{g}_t^Y}^\sharp, \omega_{\tilde{g}_t^Y}^\sharp) = n(n-1)$,  
by (\ref{ub}) we obtain that 
$$
\frac{s^{J_a}_{g^{JY}_a}}{n} = \left( {\int e^{\frac{n-2}{2}u_t}d\mu_{g_t^Y}}
\middle/ {\int e^{\frac{n+2}{2}u_t}d\mu_{g_t^Y}}\right) \mid_{t=a} < 1\, , 
$$
By continuity, there exists an $0<\varepsilon'<\varepsilon$ such that
$$
\frac{s^{J_t}_{g_t^{JY}}}{n} = \left( {\int e^{\frac{n-2}{2}u_t}d\mu_{g_t^Y}}
\middle/ {\int e^{\frac{n+2}{2}u_t}d\mu_{g_t^Y}}\right) < 1 
$$
holds for all $t \in (a-\varepsilon', a+\varepsilon')\cap [0,1]$. 

On the other hand, by the minimality of $f_{\tg_t^Y}$, and the  
first two identities in (\ref{pti}), we obtain that
$$
n(n-1)=e^{-2u_t}(n(n-1)+(n-1)(2\Delta^{\tilde{g}_t^Y}u_t - 
(n-2)\tilde{g}_t^Y(\nabla^{\tilde{g}_t^Y}u_t^\tau,\nabla^{\tilde{g}_t^Y}
 u_t^\tau))-n(n-1) \tilde{g}_t^Y(\nabla^{\tilde{g}_t^Y}u_t^\nu,
\nabla^{\tilde{g}_t^Y} u_t^\nu) 
$$
and 
$$
\| H_{f_{g_t^{JY}}}\|^2 =e^{-2u_t}n^2 \tilde{g}_t^Y(\nabla^{\tilde{g}_t^Y}
u_t^\nu, \nabla^{\tilde{g}_t^Y} u_t^\nu) \, , 
$$
respectively. By integration of the first relative to the measure 
$e^{\frac{n+2}{2}u_t} d\mu_{g_t^Y}$, we have 
$$
n(n-1)\! \int e^{\frac{n+2}{2}u_t} d\mu_{g_t^Y} =   
n(n-1) \! \int e^{\frac{n-2}{2}u_t} d\mu_{g_t^Y}-\frac{n(n-1)}{n^2}
\|H_{f_{g_t^{JY}}}\|^2 \! \int e^{\frac{n+2}{2}u_t} d\mu_{g_t^Y}\, , 
$$
from which it follows that 
$$
1+\frac{\|H_{f_{g_t^{JY}}}\|^2}{n^2} =
 \left( {\int e^{\frac{n-2}{2}u_t}d\mu_{g_t^Y}}
\middle/ {\int e^{\frac{n+2}{2}u_t}d\mu_{g_t^Y}}\right) \geq 1 
$$
for all $t$. 
\qed

\section{Sigma and almost Hermitian sigma invariants} \label{sec3} 
\subsection{The $\sigma$ invariant of the exotic spheres $M^7_k$ of
Milnor}
We begin by studying conformal properties of $Sp(2)$. We consider the 
two parameter family of left invariant metrics $g_{\lambda,\mu}$ on $Sp(2)$ that
assign norm to the elements of the Lie algebra $\mf{s}\mf{o}(2)$ by the rule 
$$
\left\| \left( \begin{array}{cr}
p & -\bar{q} \\ q & r \end{array} \right) \right\| ^2_{g_{\lambda,\mu}}=
\lambda|p|^2 + |q|^2 + \mu |r|^2\, . 
$$
The Killing form $\beta$ of Ricci tensor $\BOne$ and scalar curvature $10$ 
\cite[Corollary 7.7]{miln2} has 
$$
\left\| \left( \begin{array}{cr}
p & -\bar{q} \\ q & r \end{array} \right) \right\| ^2_{\beta}=
|p|^2 + 2|q|^2 + |r|^2\, .  
$$

The sectional curvatures $K_{\lambda,\mu}$ of $g_{\lambda,\mu}$ were computed 
in \cite[Theorem 4.3]{gracie}. In the orthonormal basis of 
$\mf{s}\mf{p}(2)$ given by  
\begin{equation} \label{basis}
\begin{array}{rrrr}
\tau_1=\left(\begin{array}{rl}
              0 & 1 \\ 
             -1 & 0 
     \end{array}\right) \, , & 
\tau_2=\left( \begin{array}{ll}
              0 & i \\ 
             i & 0 
     \end{array}\right) \, , & 
\tau_3=\left( \begin{array}{ll}
              0 & j \\ 
             j & 0 
     \end{array}\right)  \, , &  
\tau_4=\left( \begin{array}{ll}
              0 & k \\ 
             k & 0 
     \end{array}\right) \, ,  \vspace{1mm} \\
e_5=\frac{1}{\sqrt{\lambda}}\left( \begin{array}{rl}
              i & 0 \\ 
              0 & 0 
     \end{array}\right) \, , & 
e_6=\frac{1}{\sqrt{\lambda}}\left( \begin{array}{rl}
              j & 0 \\ 
              0 & 0 
     \end{array}\right) \, , & 
e_7=\frac{1}{\sqrt{\lambda}}\left( \begin{array}{rl}
              k & 0 \\ 
              0 & 0 
     \end{array}\right)  \, , &   \vspace{1mm} \\
e_8=\frac{1}{\sqrt{\mu}}\left( \begin{array}{rl}
              0 & 0 \\ 
              0 & i 
     \end{array}\right)  \, , & 
e_9=\frac{1}{\sqrt{\mu}}\left( \begin{array}{rl}
              0 & 0 \\ 
              0 & j 
     \end{array}\right)  \, , & 
e_{10}=\frac{1}{\sqrt{\mu}}\left( \begin{array}{rl}
              0 & 0 \\ 
              0 & k 
     \end{array}\right) \, ,   
\end{array}
\end{equation}
we have that
$$
\begin{array}{l}
K_{\lambda,\mu}(\tau_i,\tau_j)=4-3(\lambda+\mu)\, ,\; 1\leq i\neq j\leq 4\, ,
\vspace{1mm} \\
K_{\lambda,\mu}(e_i,\tau_j)=\lambda \, ,  \;  
K_{\lambda,\mu}(e_{i+3},\tau_j)=\mu\, ,\; i=5,6,7\, ,\; 1\leq j\leq 4\, ,
\vspace{1mm} \\
K_{\lambda,\mu}(e_i,e_j)=\frac{1}{\lambda}\, ,  \; 
5\leq i\neq j\leq 7\, , \vspace{1mm} \\
K_{\lambda,\mu}(e_i,e_j)=0\, , \; 5\leq i \leq 7\, , \; 8\leq j\leq 10\, ,
\vspace{1mm} \\
K_{\lambda,\mu}(e_i,e_j)=\frac{1}{\mu}\, , \; \text{$8\leq i \neq j \leq 10$,} 
\end{array} 
$$
respectively, and so
$$
s_{g_{\lambda,\mu}}= 2\left( \frac{3}{\lambda}+24 -6(\lambda +\mu) +
\frac{3}{\mu}\right) \, .
$$ 

The metric $g_{\frac{1}{2},\frac{1}{2}}$ on $Sp(2)$ is bi-invariant, and so
in the same conformal class as that of the Killing form $\beta$, the scale of 
the latter fixing the scale of the former. Since $Sp(2)$ acts unitarily on 
$\mb{H}^2=\mb{C}^2+\mb{C}^2=\mb{R}^2+\mb{R}^2+\mb{R}^2+\mb{R}^2$ and we have 
the consistent inclusion of projective spaces $\mb{P}^1(\mb{R}) \hookrightarrow 
\mb{P}^1(\mb{C}) \hookrightarrow \mb{P}^1(\mb{H})$ (see \cite[p. 8]{sim4}), 
$Sp(2)$ acts transitively on $\mb{S}^7$. In effect, $Sp(1)$ acts freely on 
$Sp(2)$ by the rule
\begin{equation} \label{ac0} 
q \circ \left( q_{\alpha \beta}\right) =
\left( \begin{array}{cc} q_{11} & q_{12} \\ q_{21} & q_{22} \end{array}\right)
\left( \begin{array}{cc} 1 & 0 \\ 0 & \bar{q} \end{array}\right) \, ,
\end{equation}
and the principal fibration of $Sp(2)$ induced by this action 
has base the standard $\mb{S}^7$, with the projection 
map $\pi_{\circ}$ of this fibration being given by
$$
\pi_{\circ}\left( \begin{array}{cc} q_{11} & q_{12} \\ q_{21} & q_{22} 
\end{array}\right) = \left( \begin{array}{c} q_{11} \\ q_{21}
                           \end{array}\right)
\, .
$$
The action of $Sp(2)$ on itself by left translations 
commutes with $\pi_{\circ}$, and so it is an action by
bundle morphisms of the fibration $\mb{S}^3 \hookrightarrow 
Sp(2)\stackrel
{\pi_{\circ}}{\rightarrow} {\mathbb S}^7$, thus, 
the transitive action of $Sp(2)$ on $\mb{S}^7$ is defined by matrix 
multiplication when $\mb{S}^7$ is thought of as the subset of 
$\mb{H}^2$ consisting of column vectors that have unit length, and the
isotropy group of this action is $Sp(1)$. We obtain a Riemannian 
submersion
\begin{equation} \label{rsu0}
\begin{array}{ccccc}
\mb{S}^3 & & & & \\
& \searrow & & & \\
 & & (Sp(2),\beta) & \stackrel{\pi_{\circ}}{\rightarrow} & 
({\mathbb S}^7,g) 
\end{array}
\end{equation}
with totally geodesic fibers of volume $\omega_3$. Hence, 
$$
\mu_{\beta}(Sp(2))=
\mu_g (\mb{S}^3) \times \mu_g(\mb{S}^7)=
\omega_3 \times \omega_7 = \frac{2}{3}\pi^6 \, . 
$$
Thus, $g_{\frac{1}{2},\frac{1}{2}}$ is Einstein with 
$s_{g_{\frac{1}{2},\frac{1}{2}}}=60$, and volume 
$\mu_{g_{\frac{1}{2},\frac{1}{2}}}(Sp(2))=\frac{1}{6^5}\frac{2}{3}\pi^6$.
The Yamabe invariant of its conformal class is  
$$
\lambda(Sp(2),[g_{\frac{1}{2},\frac{1}{2}}])=\lambda_{g_{\frac{1}{2},\frac{1}
{2}}}(Sp(2))=60\left( \frac{1}{6^5}\frac{2}{3}\pi^6\right)^{\frac{1}{5}} =
10 \left( \frac{2}{3}\pi^6\right)^{\frac{1}{5}}=
\lambda_{\beta}(Sp(2))\, .
$$

\begin{theorem}
If $g$ is any Riemannian metric on $Sp(2)$, then
$$
\lambda(Sp(2),[g]) \leq \lambda(Sp(2),[g_{\frac{1}{2},\frac{1}{2}}])=:
\sigma(Sp(2)) \, .
$$
\end{theorem}

{\it Proof}. For notational convenience, we let $\{ v_1, \ldots, v_{10}\}$ be 
the orthonormal basis of $\mf{s}\mf{p}(2)$ in (\ref{basis}) for the bi-invariant
metric $g_{\frac{1}{2},\frac{1}{2}}$, normalized to produce an orthonormal 
basis for the bi-invariant metric given by the ${\rm Ad}(Sp(2))$ invariant
Killing form $\beta$.
We use it to parallelize the entire manifold 
$Sp(2)$ by the left translation vector fields they generate.
If $A=(a_{ij}) \in \mb{G}\mb{L}(10,\mb{R})$, 
we define a left invariant metric $g_A$ by setting $w^A_j= 
\sum_i a_{ij}v_i$, and declaring $\{ w_1^A, \ldots, w_{10}^A\}$ to be 
an orthonormal basis for $\mf{s}\mf{p}(2)$ in the metric $g_A$. 
The scalar curvature $s_{g_A}$ is constant, and the bi-invariant metric
defined by $\beta$ is $g_{\BOne}$, and has 
$\mu_{g_{\BOne}}(Sp(2))=\frac{2}{3}\pi^6$.
We have that $\mu_{g_A}(Sp(2))=
(\lambda_1 \cdots \lambda_{10})^{\frac{1}{2}}\mu_{\BOne}(Sp(2))$, 
where $\lambda_1, \ldots, \lambda_{10}$ are the eigenvalues of $A^tA$.

If $s_{g_A}>0$, the metric $g_A$ cannot be a Yamabe metric in its class unless 
$g_A$ coincides with $g_{\BOne}$ up to homothetics. Indeed, if we assume that
$g_A$ is a Yamabe metric, we fix its volume to coincide with that of 
the Einstein metric $g_{\BOne}$, so after a scaling, $A$ is such that 
$\det{A}=1$, and 
$$ 
1=\det{(A^tA)}^{\frac{1}{2}}=(\lambda_1 \cdots \lambda_{10})^{\frac{1}{2}}\, . 
$$ 
By Theorem \ref{thM1}, the Nash isometric embeddings $f_{g_{A}}$ must 
be minimal when the metric itself is not Einstein, or scalar flat, 
in which case, the metric can be rescaled to make 
its isometric embeddings minimal also. But the optimization problem  
$$
\min{(\lambda_1\cdot \cdots \cdot \lambda_{10})} |_{\lambda_1 +\cdots +
\lambda_n=1}
$$
has solution $\lambda_i=1$, $i=1,\ldots, 10$, which leaves us at the 
extreme case of the geometric-arithmetic mean inequality. Hence, 
$f_{g_A}$ is minimal if, and only if, 
$$
1=\det{(A^tA)}^{\frac{1}{10}}=(\lambda_1 \cdots \lambda_{10})^{\frac{1}{10}}
\leq \frac{1}{10}(\lambda_1 +\cdots +\lambda_{10})=
\frac{1}{10} {\rm trace}(A^tA)= 1 \, , 
$$ 
and so $A^tA$ is similar to $\BOne$, $A\in \mb{S}\mb{O}(10)$, and
$f_{g_A}$ and $f_{g_{\BOne}}$ 
coincide up to isometries of the background sphere $(\mb{S}^{\tn},\tg)$.

Suppose now that $g^Y$ is any Yamabe metric on $Sp(2)$ of the same volume
as that of $g_{\BOne}$, and such that
\begin{equation} \label{rpn}
10=s_{g_{\BOne}} \leq  
 s_{g^Y}  < 10\cdot 9 \left( \frac{\omega_{10}}{\mu_{g_{\BOne}}}
\right)^{\frac{2}{10}}  \, . 
\end{equation}
We show that its Nash isometric embedding 
\begin{equation} \label{rpni}
f_{g^Y}: (Sp(2), g^Y) 
\rightarrow (\mb{S}^{10+p},\tg)
\hookrightarrow (\mb{S}^{\tn},\tg) 
\end{equation}
is, 
up to isometries of the background sphere, equal to $f_{g_{\BOne}}$ so 
in fact, $s_{g^Y}=s_{g_{\BOne}}$, and $\lambda(Sp(2),[g^Y])=
\lambda(Sp(2),[g_{\BOne}])$.

We can find a smooth path of volume preserving metrics $t \rightarrow g_t$ 
taking $g^Y$ at $t=0$ to $g_{\BOne}$ at $t=1$. We let $f_{g_t}$ be the path
of Nash isometric embeddings of the $g_t$s.  By Theorem    
\ref{thM2}, the path of classes $t\rightarrow [g_t]$ has a lift by Yamabe
metrics $g_t^Y$, with Nash isometric embeddings $f_{g^Y_t}$ that
deform $f_{g_t}$ by volume preserving Yamabe metrics $g_t^Y$ in 
$[g_t]$ taking $f_{g^Y}$ at 
$t=0$ to $f_{g_{\BOne}}$ at $t=1$. Then $f_{g_t^Y}$ is minimal when not 
Einstein or scalar flat, in which case then, it is minimal also after some 
scaling. Notice that $f_{g_t^Y}$ has positive scalar curvature for $t$s in a 
neighborhood $\mc{E}=(\varepsilon, 1]$ of $1$. Since $Sp(2)$ is parallelizable 
and the metric on any $f_{g_t}(Sp(2))$ is induced by that of 
$\mb{S}^{\tn}\rightarrow \mb{R}^{\tn+1}$, the tangent bundle
$T(f_{g_t}(Sp(2)))$ is an Euclidean vector bundle, and we can apply the 
Gram-Schmidt process to a fixed trivialization of this bundle, and
produce orthonormal oriented frames for the metrics $g_t^Y$ and $g_{\BOne}$,
respectively, which for $t\in \mc{E}$, can be conjugated to each other 
by $S(t) \in \mb{S}\mb{O}(10)$. Hence, up to isometries of the background,
in this range of $t$s we have that $f_{g_t^Y}= f_{g_{\BOne}}$,
$g_t^Y=g_{\BOne}$ is Einstein, and $f_{g_t^Y}$ is minimal after a 
rescaling of the metric, if necessary. By compactness of $[0,1]$, we can 
repeat this argument restricting the original path $g_t$ to say 
$[0, 1-\varepsilon/2]$, and by iteration draw this conclusion for all 
$t\in [0,1]$.
Thus, $f_{g_t^Y\mid_{t=0}}=f_{g^Y} =f_{\BOne}$ modulo isometries of the 
background sphere. 
\qed

\begin{remark}
The argument above applies for any compact Lie group $G$ with a bi-invariant 
metric $b$ of positive scalar curvature, and shows the extremality of Gromov 
of the pair $(G,b)$ in the class of Yamabe metrics \cite{grom1}. By 
the uniqueness of the solution to the Yamabe problem when $\lambda(M,[g])\leq 
0$ (cf. proof of Theorem \ref{thM2}), this can be extended to apply to $(G,b)$ 
in general. We thus derive in this manner the sigma invariant of any compact 
Lie group.
\smallskip

\begin{theorem}
Let $G$ be a compact Lie group with a bi-invariant
metric $b$. Then $\lambda(b)\geq 0$, and if $g$ is any Riemannian metric on 
$G$, we have that 
$$
\lambda(G,[g])\leq \lambda(G,[b])=\lambda_b(G))=:\sigma(G) \, .
$$  
\end{theorem}
\end{remark}
\smallskip

The spheres of Milnor \cite{mil} are the total spaces of certain 
bundles $\xi_{h\, j}$ over $\mb{S}^4$, with structure group $\mb{S}\mb{O}(4)$,
and fiber $\mb{S}^3$. If
$(h,j) \in \pi_3 (\mb{S}\mb{O}(4))=\mb{Z}+\mb{Z}$, the transition function
of $\xi_{h\, j}$ is given by the map $(h,j)\mapsto f_{hj}$, 
where $f_{hj}: \mb{S}^3 \rightarrow \mb{S}\mb{O}(4)$ is defined by 
$f_{hj}(u)v=u^h \cdot v \cdot u^j$ for $v\in {\mathbb H}$. 
If we denote by $\Sigma_{h\, j}$ the total space of the bundle $\xi_{h \, j}$, 
then if the Euler class $e(\xi_{h\, j})=h+j=1 \in H^4(\mb{S}^4,\mb{Z})$, 
$\Sigma_{h\,j}$ is homeomorphic to $\mb{S}^7$ \cite[Lemma 5]{mil}. 
For each odd integer $k$, we determine $h,j$ by the equations $h+j=1$ and 
$h-j=k$, and the Milnor seven 
dimensional sphere associated to $k$ is defined as $M^7_k=\Sigma_{h\,j}$. By 
Milnor's celebrated \cite[Theorem 3]{mil}, if 
$k^2\not \equiv 1 \mod{7}$, $M^7_k$ is homeomorphic but not diffeomorphic to 
$\mb{S}^7$. 

By taking $k=1$, we obtain as $\xi_{1\, 0}$ the Hopf fibration, and
$M^7_1=\Sigma_{1\, 0}=\mb{S}^7$. The fibration (\ref{rsu0}) induced by the 
action (\ref{ac0}) of $Sp(1)$ on $Sp(2)$ can thus be enhanced, and deformed,  
while preserving the fact that the deformed metrics on the sphere bases 
realize their sigma invariants. Indeed, we consider the family of 
left-invariants metrics 
$[0,1]\ni t \rightarrow \beta_t$ on $Sp(2)$, of
squared norm on $\mf{s}\mf{p}(2)$ defined by 
$$
\left\| \left( \begin{array}{cr}
p & -\bar{w} \\ w & r \end{array} \right) \right\| ^2_{\beta_t}=
|p|^2 + (1+t)|w|^2 + |r|^2\, .   
$$
These metrics are conformal deformations of $\beta$, of constant scalar
curvature, that yield metrics $g_{\beta_t}$ on the base $\mb{S}^7$ of the 
fibration, which in turn yield metrics $g_{\beta_t}^H$ on the base of the 
Hopf fibration. The former is given by a conformal diffeomorphism 
acting on the standard metric $g_{\beta}=g_7$ that maintains 
the geometry of 
the Hopf fiber, while the latter is given by the corresponding conformal 
diffeomorphism of $\mb{S}^4$ acting on the metric $g_4$ of the Hopf fibration 
base.
In this process, the geometry of the $\mb{S}^3$ fiber of the $Sp(2)$ fibration
remains unchanged also. By the solution to the Yamabe problem on spheres, we 
have that $g_{\beta_t}$ and $g_t^H$, whose conformal classes are those of
$[g_7]$ and $[g_4]$, respectively, are both metrics that realize the sigma 
invariants of their corresponding spheres. We thus get a path of double 
Riemannian submersions 
\begin{equation} \label{rsu0p}
\begin{array}{ccccccc}
\mb{S}^3 &  & \phantom{SSSSS} \mb{S}^3 & & & & \\ 
& \searrow & & \searrow & & &  \\
 & & (Sp(2),\beta_t)& \rightarrow & (\mb{S}^7,g_{\beta_t}) & 
\rightarrow & (\mb{S}^4, g^H_t) 
\end{array}
\end{equation}
the left end of which is a path of Riemannian submersion deformations of
(\ref{rsu0}), with metrics in its base all realizing the sigma invariant of
$\mb{S}^7$, by construction. (Notice that for $t\neq 0$, none of the 
$\beta_t$s equal any of the metrics $g_{\lambda,\mu}$.)  

By taking $k=3$, we obtain $\xi_{2\, -1}$, and
$M^7_3=\Sigma_{2\, -1}$ is an exotic sphere. 
By \cite{grme}, the action of $Sp(1)\times Sp(1)$ on $Sp(2)$ given by 
\begin{equation}
\left( \begin{array}{cc} q_1 & 0 \\ 0 & q_1 
\end{array}\right)
\left( \begin{array}{cc} q_{11} & q_{12} \\ q_{21} & q_{22} \end{array}\right)
\left( \begin{array}{cc} \bar{q}_2 & 0 \\ 0 & 1 \end{array}\right) 
\end{equation}
restricts to a free action of the diagonal $\Delta$, and 
$Sp(2)/\Delta=\Sigma_{2\, -1}=M^7_3$, the bundle $\xi_{2\, -1}$
thus induced by this action. The metrics $g_{\lambda,\mu}$ on $Sp(2)$ 
induce naturally Riemannian submersions on it.

For we let $g^{\mu}_{\lambda}$ be the metric on $\Sigma_{2\, -1}$ obtained by 
scaling the vertical space of the Hopf fibration by the factor $\lambda$. 
Then the metrics $g_{\lambda,\mu}$ descend to $g^{\mu}_{\lambda}$ on 
$\Sigma_{2\, -1}$, and
we obtain coherent Riemannian submersions 
\begin{equation} \label{rsu1}
\begin{array}{ccccccc}
\mb{S}^3 &  & \phantom{SSSSS} \mb{S}^3 & & & & \\ 
& \searrow & & \searrow & & &  \\
 & & (Sp(2),g_{\lambda,\mu})& \rightarrow & (\Sigma_{2\, -1},
g^{\mu}_{\lambda}) & \rightarrow & (\mb{S}^4, g^{H,\mu}_{\lambda}) 
\end{array}
\end{equation}
where $g^{H,\mu}_{\lambda}$ is the metric on the horizontal space of the Hopf
fibration after the $\lambda$ scaling.  
The $\mb{S}^3$ fibers are totally geodesic $3$-spheres of sectional 
curvatures $\frac{1}{\mu}$ and $\frac{1}{\lambda}$, respectively.
By construction, for $\lambda=\mu=1$, the fibers have volume
$2\pi^2$, and $\mu_{g^{H,1}_{1}}(\mb{S}^4)=\omega_4/2^4=\pi^2/6$, so 
$\mu_{g^1_{1}}(\Sigma_{2\, -1})=\pi^4/3$, and $\mu_{g_{1,1}}(Sp(2))=2\pi^6/3$.
Though the metrics $g_{\beta_0}$ and $g_1^1$ are smooth relative to distinct 
differentiable structures on the sphere, their conformal classes agree in the 
$C^0$-topology on $\mc{C}(M^7_1)$ and $\mc{C}(M^7_3)$, respectively, 
and $g_{1,1}$ and $\beta_0$ are equal volume left-invariant metrics on 
$Sp(2)$, which induce Riemannian submersions with totally geodesic fibers 
$\mb{S}^3$ of the same volume, over
equal volume seven spheres with inequivalent differentiable structures.

The entire family of Riemannian submersions $s\rightarrow (Sp(2),g_{s,s})
\rightarrow (\Sigma_{2 \, -1},g^{s}_{s})\rightarrow (\mb{S}^4,
g_{s}^{H,s})$, with geodesic fibers of equal volume and scalar
curvature $6/{\frac{1}{s}}$, is of particular interest here. Notice 
the relations
$$
s_{g_{\frac{1}{2},\frac{1}{2}}}-2\frac{6}{\frac{1}{2}} =
s_{g_{1,1}}\, , \; \; s_{g_{1,1}}+6=s_{g_7}\, , \; \; 
s_{g_{1,1}}+\frac{6}{\frac{1}{2}} =
s_{g_1^{H,1}}=s_{g_1^{H}}=48\, .
$$

\begin{theorem}
Let $M$ be any of the Milnor spheres $M^7_k=\Sigma_{h \, j}$, or any other
seven sphere with a differentiable structure. Then we have that
$$
\sigma(M)=\lambda_{g_7}(\mb{S}^7)= 42 \left( \frac{1}{3}\pi^4\right)^{\frac{2}
{7}} \, , 
$$
and the invariant is achieved by a conformal class in $\mc{C}(M)$ 
if, and only if, $M \cong \mb{S}^7$.
\end{theorem}

{\it Proof}. Since the connected sum of smooth manifolds homeomorphic to 
the sphere is itself a smooth manifold homeomorphic to the sphere, by 
\cite[\S6 p. 102]{eeku} we have that $\#_{m=1}^jM^7_3$, $1\leq j\leq 28$,
provide the $28$ distinct differentiable structures on the topological 
$7$-sphere, all of them with isomorphic combinatorial structures. By an
inductive extension of the inequality of Kobayashi in \cite[Theorem 2]{koba},
we then conclude that for any $j$ such that $1\leq j\leq 28$, we have that  
$$
\lambda_{g_7}(\mb{S}^7)\geq \sigma(\#_{m=1}^{j} M^7_3) \geq \min{\{
\sigma(M^7_3), \ldots, \sigma(M^7_3)\}} =\sigma(M^7_3) \, ,
$$
and so it suffices to prove the result just for $M=M^7_3=\Sigma_{2\, -1}$. 

We consider the family of isometric embeddings $f_{g_{\beta_t}} \in 
\mc{M}_{\mb{S}^{\tn}}(M_1^7)$ of the metrics $g_{\beta_t}$ in
(\ref{rsu0p}), and the family of isometric embeddings $f_{g_s^s}\in
\mc{M}_{\mb{S}^{\tn}}(M_3^7)$ of the metrics $g_s^s$ in (\ref{rsu1}) 
associated to the left invariant metrics $g_{s,s}$. By construction, 
in the $C^1$ topology on the space of immersions, $\lim_{t \searrow 0}f_{
\beta_{t}}=f_{g_1^1}$, while in the $C^0$ topology in the space of conformal
classes, $\lim_{t \searrow 0}[g_{\beta_{t}}]=[g_1^1]$, and we 
have that $\lambda_{g_{\beta_t}}(\mb{S}^7)=\lambda_{g_7}(\mb{S}^7)$
for all $t$, with each metric $g_{\beta_t}$ realizing the universal 
invariant upper bound. The continuity of the function 
$[g] \rightarrow \lambda(M,[g])$ follows by the continuity of
$g \rightarrow \lambda(M,[g])$. Hence,   
$\lambda(M_3^7,[g_1^1])=
\lim_{t\searrow 0}\lambda(M_1^7,[g_{\beta_t}])=
\lambda_{g_{\beta_t}}(\mb{S}^7)=\lambda_{g_7}(\mb{S}^7)$.
\qed 

\subsection{Almost Hermitian sigma invariants of products of spheres}
The sigma invariants of $M^n_k$, $1\leq k\leq n-1$, 
are all known, with 
$$
\sigma(M^n_k)=n(n-1)\omega_{n}^{\frac{2}{n}}
$$ 
if $k=1$ or $k=n-1$, but no Yamabe metric on $M^n_k$ achieves it \cite{sc2}, 
while
$$
\sigma(M^n_k)=\lambda_{g_{\mb{S}^{n,k}}}(\mb{S}^{n,k})=n(n-2)\left(\left(
\frac{k}{n}\right)^{\frac{k}{2}}\left(\frac{n-k}{n}\right)^{\frac{n-k}{2}}
\omega_k\omega_{n-k} \right)^{\frac{2}{n}}
$$
for $k$s such that $2\leq k\leq n-2$ \cite[Theorem 2]{sim4}.

\begin{theorem} \label{th12}
{\rm (a)} If $J$ is any Calabi-Eckmann complex structure on 
the product $M^n_k$ of two odd dimensional spheres of dimensions $k$ and
$n-k$, respectively, then
$$
\sigma^J(M^n_k)=\frac{1}{n-1}\sigma(M^n_k)\, , 
$$
but for no almost Hermitian Yamabe metric $g\in \mc{M}^J(M^n_k)$ is this value 
equal to $\lambda^J(g)$. 

{\rm (b)} If $J$ is an octonionic induced almost complex structure
on $M^n_k$ for $n,k$ any of the pairs $6,2$, $8,2$ or $12,6$, respectively, then
$$
\sigma^J(M^n_k)=\frac{1}{n-1}\sigma(M^n_k)=\frac{1}{n-1} 
\lambda_{g_{\mb{S}^{n,k}}}(\mb{S}^{n,k})\, ,
$$
but for no almost Hermitian Yamabe metric $g\in \mc{M}^J(M^n_k)$ is this value
equal to $\lambda^J(g)$. 
\end{theorem}

{\it Proof}. By \cite[\S2]{sc2}, there exists a sequence $a \rightarrow r_a 
\nearrow 1$, and associated sequence of Yamabe metrics $g_a$ in the 
conformal class of the product metric $g^n_1(r_a)$ on the Calabi-Eckmann 
manifold $M_1^{n}(r_a,\sqrt{1-r_a^2}) \hookrightarrow \mb{S}^{n+1}$, such that
the sequence of Yamabe conformal invariants of $g_a$  
increases to the universal bound $n(n-1)\omega_{n}^{\frac{2}{n}}$, which
is not realized by any Yamabe metric on the manifold. Since all of these
manifolds are locally conformally flat, by Theorem \ref{th5}, the metrics 
$g_a$ are $J_{n,1}^{r_a}$ Yamabe metrics also, and we have that  
$\lambda^{J_{n,1}^{r_a}}(g_a))=\lambda(g_a)/(n-1) \nearrow
n\omega_{n}^{\frac{2}{n}}$, but the limit bound is not
realized by any almost Hermitian Yamabe metric. This
proves (a) for $k=1$, and by interchanging the sphere factors, for $k=n-1$
also.

The remaining cases in (a), and all of the three cases in (b), for which we 
have $n\geq 6$,
follow by the same argument. Notice that in general, by (\ref{rem}), we have 
the inequalities
\begin{equation} \label{iubjy}
(n-1) \lambda^J(M,[g]) \geq \lambda(M,[g]) + 2(n-1)\inf_{g\in \mc{M}^J_{[g]}(M)}
\frac{1}{\mu_{g}^{\frac{2}{N}}}\int W_{g}(\omega^{\sharp}_{g},
\omega^{\sharp}_{g})d\mu_g \, ,  
\end{equation}
and 
\begin{equation} \label{iuby}
\lambda(M,[g]) \geq (n-1)\lambda^J(M,[g]) +2(n-1)\inf_{g\in \mc{M}^J_{[g]}(M)} 
 \frac{-1}{\mu_{g}^{\frac{2}{N}}}\int W_{g}(\omega^{\sharp}_{g},
\omega^{\sharp}_{g}) d\mu_g \, ,  
\end{equation}
respectively, and so by the trychotomy of cases that solve the almost
Hermitian Yamabe problem \cite{rss}, the classes $[g]\in \mc{C}^{J}(M^n_k)$ 
with the largest values of $\lambda^J(M^n_k,[g])$ are found among those for
which $\lambda(M^n_k,[g])>0$, and which admit a representative $g$ such that 
$W_{g}(\omega^{\sharp}_{g},\omega^{\sharp}_{g})\geq 0$, and is positive at some
point.

Thus, let us consider $g\in \mc{M}^J(M^n_k)$ such that $\lambda(M^n_k,[g])>0$
with $W_g(\omega_g^{\sharp},\omega_g^{\sharp})\geq 0$ nontrivial.
We let $g^Y$ and $g^{JY}$ be Yamabe and almost Hermitian Yamabe metrics in 
$[g]$ of volume $\mu_g$ each. By Theorem \ref{thM3}, there exists a 
function $u$ defined on a tubular neighborhood of $f_{g^Y}(M^n_k)
\hookrightarrow \mb{S}^{\tn}$ such that on points of $f_{g^Y}(M^n_k)$,
we have that
$$
s^J_{g^{JY}}= e^{-2u}\left( \left( \frac{1}{n-1}s_{g^Y} 
+ 2\Delta^{g^Y} u -(n-2)\tg(du,du) \right)
+2 W_{g^Y}(\omega^{\sharp}_{g^Y}, \omega^{\sharp}_{g^Y})  
\right)\, , 
$$
whose intrinsic meaning on $M^n_k$ provided with the metric
$g^{JY}=e^{2u}g^Y$ is just that
$$
(n-1)s^J_{g^{JY}}=s_{e^{2u}g^Y}+ 2(n-1)e^{-2u}
W_{g^Y}(\omega^{\sharp}_{g^Y}, \omega^{\sharp}_{g^Y}) \, . 
$$
The nonnegative function $W_{g^Y}(\omega^{\sharp}_{g^Y}, 
\omega^{\sharp}_{g^Y})$ is nontrivial, and we have that
$$
W_{g^{JY}}(\omega^{\sharp}_{g^{JY}}, \omega^{\sharp}_{g^{JY}}) =
e^{-2u}
W_{g^Y}(\omega^{\sharp}_{g^Y}, \omega^{\sharp}_{g^Y}) \, . 
$$
We prove that
\begin{equation} \label{weq}
\inf_{g\in \mc{M}_{[e^{2u}g^Y]}}
\frac{1}{\mu_{g}^{\frac{2}{N}}}\int W_{g}(\omega^{\sharp}_{g},
\omega^{\sharp}_{g})d\mu_{g}=0 \, .  
\end{equation}
We proceed along the lines of \cite[\S6 \& 7]{lepa}, in turn inspired by
the crucial work of Schoen \cite{sc} (and Aubin \cite{au}), to make suitable 
choices of test functions to accomplish this.  

Since $\lambda(M^n_k,[g^Y])>0$, the Green function $\Gamma_p$ at each 
$p\in M^n_k$ for the conformal Laplacian exists and is positive. We choose 
$p$ to be a point where $W_{g^Y}(\omega^{\sharp}_{g^Y}, \omega^{\sharp}_{g^Y})$ 
achieves its maximum, and let $(\hat{M},\hat{g})$ be the stereographic 
projection of $M^n_k$ from $p$ provided with the asymptotically flat metric 
$\hat{g}=G^{\frac{4}{n-2}}g^Y$, $G=4\omega_n\Gamma_p$ (see \cite[p. 64]{lepa}).
Working in inverted normal coordinates $x$, we set the radial distance
function $r(x)=|x|$ as a (extended) smooth positive funtion on $\hat{M}$, and
use it to consider
$$ 
u_\alpha(x)=\left(\frac{|x|^2+\alpha^2}{\alpha}\right)^{(2-n)/2}\, .
$$ 
If $\hat{M}$ were equal to $\mb{R}^n$, these functions would all 
realize the optimal Sobolev constant 
$$ 
a \|\nabla u_\alpha\|^2_{L^2({\mathbb R}^n)}=\Lambda 
\|u_\alpha\|^2_{L^N({\mathbb R}^n)}\, , 
$$ 
where $a=4(n-1)/(n-2)$ and $\Lambda=n(n-1)\omega_n^{2/n}$. Since the
metric distortion coefficient of the asymptotically flat metric $\hat{g}$ is 
positive, the $L^N$ norm of $u_{\alpha}$ in the metrics $\hat{g}$ can be
estimated from above and below by the $L^N$ norm of $u_{\alpha}$ in the 
metric $g^Y$. If $g_{\alpha}=u_{\alpha}^{N-2}g^Y$, by the conformal properties 
of the Weyl tensor, we have that
$$
\int_M W_{g_{\alpha}}(\omega^{\sharp}_{g_{\alpha}},\omega^{\sharp}_{g_{\alpha}})
d\mu_{g_{\alpha}} =
\int_M u_{\alpha}^{2}W_{g^Y}(\omega^{\sharp}_{g^Y},\omega^{\sharp}_{g^Y})
d\mu_{g^Y} \leq W(p) \int_M u_{\alpha}^{2}
d\mu_{g^Y} \, ,
$$
and since $s_{g^Y}>0$, so the volume of geodesic spheres grow slower than in
$\mb{R}^n$, we conclude that
$$
\int_M W_{g_{\alpha}}(\omega^{\sharp}_{g_{\alpha}},\omega^{\sharp}_{g_{\alpha}})
d\mu_{g_{\alpha}} 
\leq W(p) \int_{\mb{R}^n} u_{\alpha}^{2}(x) dx 
\leq W(p)\omega_{n-1}\alpha^2 I_n \, ,
$$
where 
$$
I_n=\int_{0}^{\infty} s^{n-1}(1+s^2)^{2-n}ds = \frac{2^{3-n}\pi^{\frac{1}{2}}}{(
n-4)} \frac{\Gamma\left( \frac{n}{2}\right)}
{\Gamma\left( \frac{n-1}{2}\right)} \, .
$$
Thus, for all sufficiently small $\alpha$s, there exist a constant $C$ such 
that
$$
\frac{1}{\mu_{g_\alpha}^{\frac{2}{N}}}\int_M W_{g_{\alpha}}
(\omega^{\sharp}_{g_{\alpha}},\omega^{\sharp}_{g_{\alpha}})
\leq 
\frac{C}{\mu_{g_\alpha}^{\frac{2}{N}}}W(p)\omega_{n-1}\alpha^2 \, , 
$$
and as the metrics $e^{2u}g^Y$ and $g^Y$ are in the same conformal class,  
(\ref{weq}) follows by choosing $\alpha \searrow 0$.

By (\ref{iubjy}), (\ref{weq}) implies that
$(n-1)\lambda^{J}(M^n_k,[g]) \geq \lambda(M^n_k,[g]) > 0$, 
with the conformal invariants achieved at different metrics 
$g^{JY}$ and $g^Y$. On the other hand, by (\ref{iuby}) and the nonnegativity
of $e^{-2u}W_{g^Y}(\omega^{\sharp}_{g^Y}, \omega^{\sharp}_{g^Y})$, we have that
$\lambda(M^n_k,[g]) \geq (n-1)\lambda^J(M^n_k,[g])$. Therefore, 
$$
\lambda(M^n_k,[g]) = (n-1)\lambda^J(M^n_k,[g]) \, .
$$
Since the sigma invariant
of $M^n_k$ is realized by $g_{\mb{S}^{n,k}} \in \mc{M}(M^n_k)\cap 
\mc{M}^J(M^n_k)$, the 
supremum of $\lambda(M^n_k,[g])$ over $\mc{C}(M^n_k)$ is equal to
the supremum of $\lambda(M^n_k,[g])$ over $\mc{C}^J(M^n_k)$,
and we obtain that 
$$
(n-1)\sigma^J(M^n_k)=(n-1)\sup_{[g]\in \mc{C}^J}\lambda^J(M^n_k,[g]) =
\sup_{[g]\in \mc{C}}\lambda(M^n_k,[g])= \lambda(g_{\mb{S}^{n,k}})=\sigma
(M^n_k)\, .
$$
By Theorem \ref{th5}, the almost Hermitian invariant $\sigma^J(M^n_k)$ cannot
be realized by any almost Hermitian Yamabe metric $g$ in $\mc{M}^J(M^n_k)$. 
\qed

\begin{remark} \label{re13}
The argument above for (b) fails when $n=4$ and the manifold is the
K\"ahler product $M^4_2$ of two spheres. This is related to the
structure of $(\mb{S}^{4,2},g_{\mb{S}^{4,2}})$ as a symplectic ruled surface  
\cite[Theorem \S 2.4.C]{grom}, a manifold that carries connected regular
$J$-holomorphic curves of genus $1$ representing the homology class 
$p[\mb{S}^2]+ q[\mb{S}^2]\in H_2(M^4_2,\mb{Z})$, $p,q$ any pair of 
nonnegative integers, and whose group of symplectomorphisms in the identity component of the
diffeomorphism group is connected, and deformation retracts to $\mb{S}\mb{O}(3)
\times \mb{S}\mb{O}(3)$. But for product metrics $g^4_2(r)$ on $M^4_2(r,\sqrt{1
-r^2})$ with $r^2\neq \frac{1}{2}$, Gromov indicates that
the analogous group acquires a new infinite
order element in its fundamental group (see \cite[Proposition 1.2]{mmcd}), and 
so we may produce a family $J_t$ of almost complex structures
associated to the tori representatives of the homology classes  
$[\mb{S}^2]+ 2[\mb{S}^2]$, or $2[\mb{S}^2]+ [\mb{S}^2]$, 
and $J_t$ compatible metrics $g_t$, with symplectic forms $\omega_t$ in the 
real cohomology classes $[\mb{S}^2]+ 2t[\mb{S}^2]$, or 
$2t[\mb{S}^2]+ [\mb{S}^2]$, respectively, such that 
$s_{g_t}^{J_t}\rightarrow 12$ and 
$\mu_{g_t}(M^4_2)\rightarrow \omega_4=\frac{8}{3}\pi^2$ as $t\rightarrow 1$.
This suggests that
$\sigma^J(M^4_2)=4\omega_4^{\frac{1}{2}}=8\pi \left(\frac{2}{3}
\right)^{\frac{1}{2}}$, without this value being realized by an almost 
Hermitian Yamabe metric. As proved earlier in (a), this is also the value of 
the almost Hermitian sigma invariant of $M^4_1$, a manifold that is the total 
space of Riemannian submersions over $\mb{S}^2=\mb{P}^1(\mb{C})$, with scalar 
flat tori fibers.
\end{remark}

\subsection{The conformal Pascal triangle of Yamabe manifolds}
It is natural to extend the definition of the Yamabe functional of a metric on
$M^n$ as being identically zero when $n=0$ or $1$, or 
$4\pi \chi(M)/(\frac{1}{4}\mc{W}(M,[g])^{\frac{1}{2}}$, 
where $[g]$ is the conformal class of $M$ of the
smallest Willmore energy $\mc{W}(M,[g])$ when $n=2$,  
respectively. In the former cases, the manifolds have trivial sigma invariants,
realized by the oriented point deformation 
retract of a totally geodesic oriented $2$-disk of radius $\pi/2$ in 
$\mb{S}^{\tn}$ and its bounding geodesic circle, respectively, 
while on surfaces, the Yamabe invariant of 
a conformal class is given by
$4\pi \chi(M)/(\frac{1}{4}\mc{W}(M,[g])^{\frac{1}{2}}$,  
and in a precise sense, the sigma invariant of $M$ is 
achieved by the conformal class with the most symmetries in its fundamental 
domain \cite[Theorems 1,8,9]{sim2}, for instance, 
$\sigma(\mb{S}^1\times \mb{S}^1)=0  < \sigma(\mb{P}^2(\mb{R}))=
\frac{2\sqrt{6}}{3}\pi^{\frac{1}{2}} < 
\sigma(\mb{S}^2)=4\pi^{\frac{1}{2}}$, with the Clifford torus 
$(\mb{S}^{2,1},g_{S^{2,1}})$ realizing the first of these three terms. 

Yamabe Riemannian manifolds $(M,g)$ then fit into a conformal Pascal triangle 
arranged according to the symmetries of $g$. We 
sketch its description by using the Yamabe manifolds $(\mb{S}^n,g)$ and 
$(\mb{S}^{n,k},g_{\mb{S}^{n,k}})$ as its initial elements. 
For we consider the arrangement of these manifolds given by  
$$
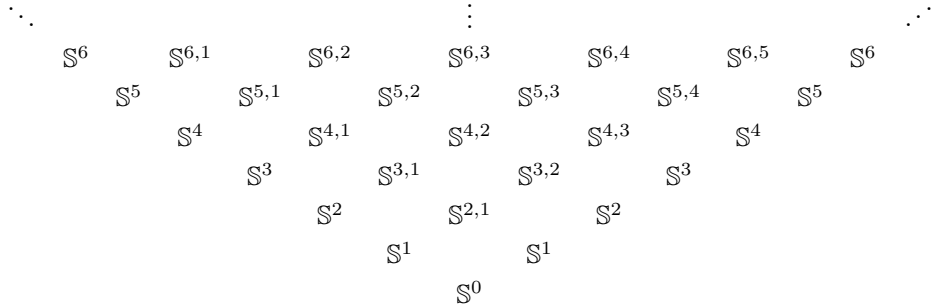
\begin{array}{ccccccccccccccccc}
 & \ddots & & & & & & & \vdots & & & & & &
 & \iddots \vspace{1mm}\\
 & & \mb{S}^6 & & \mb{S}^{6,1} & & \mb{S}^{6,2} & & \mb{S}^{6,3} & & 
 \mb{S}^{6,4} & &  \mb{S}^{6,5} & & \mb{S}^6 & \vspace{1mm} \\
 & & & \mb{S}^5 & & \mb{S}^{5,1} & & \mb{S}^{5,2} & & \mb{S}^{5,3} & & 
 \mb{S}^{5,4} & &  \mb{S}^5 & & & \vspace{1mm} \\
 & & & & \mb{S}^4 & & \mb{S}^{4,1} & & \mb{S}^{4,2} & & \mb{S}^{4,3} & & 
\mb{S}^4 & & & &  \vspace{1mm} \\
 & & & & & \mb{S}^3 &  &  \mb{S}^{3,1} & & \mb{S}^{3,2} & & \mb{S}^3 & & & & & 
 \vspace{1mm} \\
 & & & & & & \mb{S}^2 & & \mb{S}^{2,1} & & \mb{S}^2 & & & & & &  \vspace{1mm} \\
 & & & & & & & \mb{S}^1 & & \mb{S}^1 & & & & & & &   \vspace{1mm} \\
 & & & & & & & & \mb{S}^0 & & & & & & & & 
\end{array}
$$
The conformal Pascal triangle corresponds to this arrangement, where
each entry $M$, on each row, is replaced by the Yamabe functional value of 
the Yamabe metric it carries. A Yamabe manifold $(M^n,g)$ is hanged up off the 
$n$th horizontal line in the triangle, at a location determined by its 
symmetries, and at a height equal to $\lambda(g)=\lambda(M^n,[g])$. 
All the manifolds depicted explicitly in the $n$th line of the triangle above,
with one exception, are those appearing at both ends of Simons' gap theorem 
\cite[Theorem 5.3.2, Corollary 5.3.2]{simo}, \cite[Main Theorem]{cdck}, 
\cite[Corollary 2]{blaine}, that is to say, manifolds $(M^n,g)$ with minimal 
isometric embeddings into $\mb{S}^{n+p}$ such that 
$\| \alpha \|^2 \leq np/(2p-1)$, 
and for which the theorem implies that either $\| \alpha \|^2 =0$, in 
which case $p=1$, and $M=\mb{S}^n$, the case at the extreme edges of the 
triangle, or $\| \alpha \|^2=n$, and either $p=1$, and $M$ is any of 
$\mb{S}^{n,k}$, $k=1, \ldots, n-1$, the cases within the said extreme edges 
depicted, or $p=2$, $n=2$, and $M=\mb{P}^2(\mb{R})$ with a metric of scalar 
curvature $2/3$, and area $6\pi$, isometrically embedded into $\mb{S}^4$ 
(the non depicted exception). For $k$s in the range 
$1\leq k \leq \left[ \frac{n}{2}\right]$, the sequence 
$\lambda(g_{\mb{S}^{n,k}})$ is decreasing \cite[Example 10(1)]{sim5}, and 
by Aubin's universal bound, each of its terms is bounded above by $\lambda(g)=
\lambda(\mb{S}^n,[g])$. The isometry given by the interchange of the factor
spheres in the product makes of $\mb{S}^{n,n-k}$ a symmetric mirror image 
of $\mb{S}^{n,k}$, which hangs on the mirror side of the first half cut of the
triangle through its vertical middle axis, at the same height.

The sign of the hanging level of a manifold with a metric realizing its 
sigma invariant is determined by the trychotomy theorem of Kazdan \& Warner 
\cite{kawa}: The hanging level of a manifold of type I, with a Yamabe metric 
that achieves its sigma invariant, is positive; the hanging level of manifolds 
of type II, with scalar flat metrics on them, is zero; and the hanging level of 
manifolds of type III, with any Yamabe metric on them, is negative. Notice the
implications of these facts for manifolds whose sigma invariant is not 
achieved by any conformal class, whose statement we bypass.

The hanging level for known manifolds of positive sigma invariant of the 
same dimension follows the observed pattern in the extended definition of the 
concept for manifolds of dimension two (or less) above. For instance, 
by \cite[Theorems 2 \& 5]{sim4}, we have that
$$
\begin{array}{ccccccc}
\sigma(\mb{P}^4(\mb{R})) & < & \sigma(\mb{S}^{4,2})  & < &  
\sigma(\mb{P}^2(\mb{C})) & < & 
\sigma(\mb{P}^1(\mb{H})) = \sigma(\mb{S}^4)\, , \\
\sigma(\mb{P}^6(\mb{R})) & < & \sigma(\mb{P}^3(\mb{C})) & < & 
\sigma(\mb{S}^{6,3}) & < & \sigma(\mb{S}^6)\, , \\
\sigma(\mb{P}^8(\mb{R})) & < & \sigma(\mb{P}^4(\mb{C})) & < & 
\sigma(\mb{P}^2(\mb{H})) & < &  
\sigma(\mb{S}^{8,4}) \hspace{2mm}  <  
\hspace{2mm} \sigma(\mb{S}^8)\, .
\end{array}
$$

By Theorem \ref{th12}, we can formulate an analogous 
almost Hermitian version of this conformal Pascal triangle for 
almost Hermitian Yamabe manifolds, the edges of this version hanging at height
given by the limit of $(n-1)\lambda^{J_{n,1}}(g^n_1(r_a))$ for the 
Calabi-Eckmann 
manifolds $(M^n_1(r_a,\sqrt{1-r_a^2}),g^n_1(r_a))$ in the early part 
of the proof of the theorem, which is equal to the universal
upper bound $\sigma(\mb{S}^n)$ of the hanging level of the edges in the one 
above.  It is remarkable that
if the qualitative argument of Remark \ref{re13} can be formalized, the 
K\"ahler surface $(\mb{S}^{4,2}, g_{\mb{S}^{4,2}})$ could be 
deformed to produce sequences paralleling these 
$(M^4_1(r_a,\sqrt{1-r_a^2}),g^4_1(r_a))$s, or their symmetric images, whose
almost Hermitian Yamabe invariants times $3$ approach the hanging level of 
the edges of the Pascal triangle for Yamabe manifolds of dimension $n=4$, 
the factor that gets closer to a representative of the homology class
$2[\mb{S}^2]$ determining the side edge of the triangle approached by each of 
these sequences. There is no 
counterpart of this for manifolds of dimension $n=2$ in the triangle 
because the Clifford torus $\mb{S}^{2,1}$ is flat, and hangs at height $0$,
as does any torus in a different conformal class with its standard minimally 
embedded flat metric, while the edge manifolds $(\mb{S}^2,g)$, though satisfying
the $3$ geodesics theorem in abundance, hang both at the same positive height
$4\pi^{\frac{1}{2}}$. The two components of the initial 
$\mb{S}^{0}$ in the triangular arrangement, correspond to the two generators of 
$\pi_{1}(\mb{S}^1,*)=H_1(\mb{S}^1,\mb{Z})$ of its two circle edge elements at 
the $n=1$ horizontal level. On the other hand, in real dimension $n=4$, 
we have that
$$
\sigma(\mb{S}^{4})= 
\sigma(\mb{P}^1(\mb{H})) > 
\sigma(\mb{P}^{2}(\mb{C})) > \sigma(\mb{S}^{4,2}) >
\sigma(\mb{P}^4(\mb{R}))\, ,
$$ 
while 
$$
\sigma(\mb{S}^{n})>  
\sigma(\mb{S}^{n,\frac{n}{2}})> \sigma(\mb{P}^\frac{n}{2}(\mb{C})) > 
\sigma(\mb{P}^n(\mb{R}))\, , \; \text{for $n=2k$, $k=3,4,\ldots$,}
$$ 
and 
$$
\sigma(\mb{S}^{n})>  
\sigma(\mb{S}^{n,\frac{n}{2}})> 
\sigma(\mb{P}^{\frac{n}{4}}(\mb{H})) > \sigma(\mb{P}^{\frac{n}{2}}(\mb{C})) > 
\sigma(\mb{P}^{n}(\mb{R}))
\, , \; \text{for $n=4k$, $k=2,3,\ldots$,}
$$
respectively, in an apparent reversal of the location in the ordering of 
the product $\mb{S}^{n,\frac{n}{2}}$ 
pass the $n=4$ case, but which can be explained by the CW decompositions of 
the manifolds involved, and how their cells, with the induced metrics on them, 
fit into one another by natural inclusions, generating cycles of 
{\it canonical shape} (cf. \cite[\S 4]{gracie}) . The low dimensional 
topology phenomena, and the rest of the patterns in higher dimension, are 
suited for explanations in terms of the rational homotopy of the manifolds 
involved, which precedes torsion considerations, but goes well beyond 
their de Rham cohomology rings \cite{sull1,sull}.

\end{document}